\documentclass[aos]{imsart}

\RequirePackage{amsthm,amsmath,amsfonts,amssymb}
\RequirePackage[authoryear]{natbib} %% uncomment this for author-year bibliography
\RequirePackage[colorlinks,citecolor=blue,urlcolor=blue]{hyperref}
\RequirePackage{graphicx}

\startlocaldefs
\newtheorem{theorem}{Theorem}[section]
\newtheorem{lemma}[theorem]{Lemma}
\newtheorem{proposition}[theorem]{Proposition}

\endlocaldefs

\usepackage{graphicx} % more modern
\usepackage{subfigure} 
\usepackage{placeins} %needed for floatbarrier to control mix of figs and biblio at end
\usepackage{natbib}
\usepackage{comment}
\usepackage{enumerate}
\usepackage{dblfloatfix}
\usepackage{hyperref}

\usepackage[symbol]{footmisc}

\usepackage{mathtools}

\usepackage{amsmath}
\usepackage{amssymb}
\usepackage{amsthm}
\usepackage{upgreek}
\usepackage{mathrsfs}

\newcommand{\R}{\mathbb{R}}

\renewcommand{\hat}[1]{\widehat{#1}}

\renewcommand{\P}{\mathbb{P}}
\newcommand{\E}{\mathbb{E}}

\newcommand{\e}{\epsilon}
\newcommand{\ve}{\varepsilon}

\newcommand{\var}{\operatorname{var}}

\newcommand{\ttop}{^{\top}}

\newcommand{\ts}{\textstyle}

\newcommand{\Log}{\textup{log}}

\newcommand{\mnorm}[1]{\left\vert\kern-1.5pt\left\vert\kern-1.5pt\left\vert #1\right\vert\kern-1.5pt\right\vert\kern-1.5pt\right\vert}

\def\references{\bibliography{rootnCLTbib_2nd_revision.bib}}

\begin{document}

\begin{frontmatter}
\title{ Central Limit Theorem and Bootstrap Approximation in High Dimensions: Near  \lowercase{$1/\sqrt{n}$} Rates via Implicit Smoothing}

\runtitle{Near \lowercase{$1/\sqrt{n}$} Rates for CLT and Bootstrap in High Dimensions}

\begin{aug}
%%%%%%%%%%%%%%%%%%%%%%%%%%%%%%%%%%%%%%%%%%%%%%
%%Only one address is permitted per author. %%
%%Only division, organization and e-mail is %%
%%included in the address.                  %%
%%Additional information can be included in %%
%%the Acknowledgments section if necessary. %%
%%%%%%%%%%%%%%%%%%%%%%%%%%%%%%%%%%%%%%%%%%%%%%
\author[A]{\fnms{Miles E.} \snm{Lopes}\thanksref{t1}\ead[label=e1]{melopes@ucdavis.edu}},
%\author[B]{\fnms{Second} \snm{Author}\ead[label=e2,mark]{second@somewhere.com}}
%\and
%\author[B]{\fnms{Third} \snm{Author}\ead[label=e3,mark]{third@somewhere.com}}
%%%%%%%%%%%%%%%%%%%%%%%%%%%%%%%%%%%%%%%%%%%%%%
%% Addresses                                %%
%%%%%%%%%%%%%%%%%%%%%%%%%%%%%%%%%%%%%%%%%%%%%%
\address[A]{University of California, Davis,
\printead{e1}}

%\address[B]{Department,
%University or Company Name,
%\printead{e2,e3}}
\end{aug}

\begin{abstract}
Non-asymptotic bounds for Gaussian and bootstrap approximation have recently attracted significant interest in high-dimensional statistics.
This paper studies Berry-Esseen bounds for such approximations with respect to the multivariate Kolmogorov distance, in the context of a sum of $n$ random vectors that are $p$-dimensional and i.i.d. Up to now, a growing line of work has established bounds with mild logarithmic dependence on $p$. However, the problem of developing corresponding bounds with near $n^{-1/2}$ dependence on $n$ has remained largely unresolved. Within the setting of random vectors that have sub-Gaussian or sub-exponential entries, this paper establishes bounds with near $n^{-1/2}$ dependence, for both Gaussian and bootstrap approximation. In addition, the proofs are considerably distinct from other recent approaches and make use of an ``implicit smoothing'' operation in the Lindeberg interpolation.

\end{abstract}

\thankstext{t1}{Supported in part by NSF grant DMS-1915786.}



\end{frontmatter}
%%%%%%%%%%%%%%%%%%%%%%%%%%%%%%%%%%%%%%%%%%%%%%
%% Please use \tableofcontents for articles %%
%% with 50 pages and more                   %%
%%%%%%%%%%%%%%%%%%%%%%%%%%%%%%%%%%%%%%%%%%%%%%
%\tableofcontents

\section{Introduction}\label{sec:intro}
The analysis of Berry-Esseen bounds for Gaussian and bootstrap approximation has become a quickly growing topic in high-dimensional statistics. Indeed, much of the work in this direction has been propelled by the fact that such approximations are essential tools for a wide variety inference problems.

To briefly review the modern literature on multivariate Berry-Esseen bounds, a natural starting point is the seminal paper~\citep{Bentkus:2003}. In that work, Bentkus studied Gaussian approximation of a sum $S_n=n^{-1/2}\sum_{i=1}^n X_i$ of centered i.i.d.~random vectors in $\R^p$ with an identity covariance matrix. Letting $Y$ denote a centered Gaussian random vector with $\E[YY\ttop]=\E[X_1X_1\ttop]$, and letting $\mathscr{A}$ denote the class of all Borel convex subsets of $\R^p$, Bentkus' work showed that under suitable moment conditions, the distance \smash{$\sup_{A\in\mathscr{A}}|\P(S_n\in A)-\P(Y\in A)|$} is at most of order $p^{7/4}n^{-1/2}$. (See also~\citep{Bentkus:2005,Raic:2019} for refinements and further references, as well as the recent work~\citep{Fang:2020:convex}, which improved the rate to $p^{5/4}n^{-1/2}$, up to a logarithmic factor.) However, despite the strength of such results, they are not directly applicable to situations where $p$ is larger than $n$.

In high-dimensional settings, the paper~\citep{CCK:2013} achieved a breakthrough by demonstrating that if $\mathscr{A}$ is taken instead to be a certain class of hyperectangles, then the corresponding distance can be bounded at a rate that has a logarithmic dependence on $p$, such as $\log^{7/8}(pn)n^{-1/8}$, and similarly for bootstrap approximation. Subsequently, the  papers~\citep{CCK:2017} and~\citep{CCKK:2019}, showed that when $\mathscr{A}$ includes all hyperrectangles, the rates for Gaussian and bootstrap approximation can be improved to $\log^{7/6}(pn)n^{-1/6}$
 and $\log^{5/4}(pn)n^{-1/4}$ respectively. Meanwhile, a parallel series of works ~(\cite{Deng:Zhang:2020, Kuchibhotla:2018, Koike:2019, Deng:2020, Das:Lahiri:2020}) 
has made strides in showing that Gaussian and bootstrap approximation can succeed asymptotically when $\log^{\kappa}(p)=o(n)$ and $2\leq \kappa\leq 5$.

With regard to possible improvements beyond near $n^{-1/4}$ rates, it is well known from classical results that the $n^{-1/2}$ rate in the univariate Berry-Esseen theorem is the best that can be hoped for in general. (More formally, there is a distribution for a scalar variable $X_1$ such that the quantity $\sup_{t\in\R}n^{1/2}|\P(S_n\leq t)-\P(Y\leq t)|$ is bounded away from 0 as $n\to\infty$~\citep{Esseen:1956}.)
For this reason, there has been significant interest in finding out whether near $n^{-1/2}$ rates of Gaussian approximation are achievable in high-dimensional situations. Two papers that have recently made progress on this question are~\citep{Lopes:2020} and~\citep{Fang:Koike:2020}.
The first of these papers considered a setting of ``weak variance decay'', where $\var(X_{1j})=\mathcal{O}(j^{-a})$ for all $1\leq j\leq p$, with $a>0$ being an arbitrarily small parameter. Under this type of structure, the authors
established the rate $n^{-1/2+\delta}$ for arbitrarily small $\delta>0$, when $\mathscr{A}$ is a certain class of hyperrectangles. In a different direction, the paper \citep{Fang:Koike:2020} dealt with a setting where $\mathscr{A}$ includes all hyperrectangles, establishing the rate $\log^{3/2}(p)\log(n)n^{-1/2}$
 when $X_1$ has a log-concave density, as well as the rate $\log^{4/3}(pn)n^{-1/3}$ when $X_1$ has sub-Gaussian entries (but need not have a density).
In the current work, we focus on the latter case, and the main contribution of our first result (Theorem~\ref{thm:main}) is in establishing a rate with near $n^{-1/2}$ dependence on $n$.

In addition to the work on Gaussian approximation described above, there are a few special cases where near $n^{-1/2}$ rates are known to be achievable via bootstrap approximation. First, in the setting of weak variance decay, it was shown in~\citep{Lopes:2020} that the mentioned $n^{-1/2+\delta}$ rate holds for bootstrap approximation as well. (See also~\citep{Lopes:2022} for more recent extensions of such results.) Second, the paper~\citep{CCKK:2019} showed that near $n^{-1/2}$ rates can be achieved when bootstrap methods are used in particular ways. Namely, this was demonstrated in the case when the data have a symmetric distribution and Rademacher weights are chosen for the multiplier bootstrap, or when bootstrap quantiles are adjusted in a conservative manner. (See also~\citep{Deng:2020} for further work in this direction.) In relation to these results, the current paper makes a second contribution in Theorems~\ref{thm:boot} and~\ref{thm:bootemp} by showing that near $n^{-1/2}$ rates of bootstrap approximation hold without variance decay, symmetry, or conservative adjustments.

Perhaps the most important point to discuss concerning the proofs is the use of smoothing techniques. As is well known, these techniques are based on using a smooth function, say $\psi:\R^p\to\R$, depending on a set $A\subset\R^p$, such that \smash{$\E[\psi(S_n)]\approx \P(S_n\in A)$.} Although these techniques are of fundamental importance, one of their  drawbacks is that they often incur an extra smoothing error $|\P(S_n\in A)-\E[\psi(S_n)]|$, which must be balanced with errors from various other approximations. Moreover, this balancing process often turns out to be a bottleneck for the overall rate of distributional approximation. 

As a way of avoiding this bottleneck, we use a smoothing function that arises ``implicitly'' as part of the Lindeberg interpolation scheme---which has the benefit that it does not create any smoothing error. More concretely, if $X_1,\dots,X_n$ are non-Gaussian and if $Y_1,\dots,Y_n$ are Gaussian, then this notion of smoothing is based on the fact that the probability 
$\P(\sum_{i=1}^{k}X_i+\sum_{j=k+1}^nY_j\in A)$
 can be equivalently written as $\E[\tilde\psi(\sum_{i=1}^k X_i)]$, for a particular smooth random function $\tilde\psi$ defined in terms of $Y_{k+1},\dots,Y_n$. (See Section~\ref{sec:smooth} for further details.) To a certain extent, this approach is akin to the ``method of compositions'' \citep[cf.][]{Senatov:1981}, but it is not exactly the same, since the latter approach still creates a source of smoothing error. Another notable feature of the current approach is that the derivatives of $\tilde\psi$ 
can be controlled effectively, due to the results in~\citep{Bentkus:1990,Anderson:1998,Fang:Koike:2020}. However, by itself, the use of implicit smoothing does not seem to provide a way to handle all aspects of the Lindeberg interpolation, 
  and a second important ingredient in the proof is the use of induction.
  In particular, the use of induction here is influenced by the paper~\citep{Bentkus:2003}, even though the approach to smoothing in that work is different.\\

\noindent\emph{Remark.} After the initial version of this work appeared in~\citep{Lopes:2020CLT}, the results were extended in~\citep{Kuchibhotla:2020,CCK:2020}. Improved aspects of the results include relaxed tail assumptions, and refined logarithmic factors in the rates of approximation.\\

\noindent\emph{Notation.} For a random variable $U$ and a number $q\in\{1,2\}$, define the $\psi_q$-Orlicz norm $\|U\|_{\psi_q}=\inf\{t>0| \E[\exp(|U|^q/t^q)]\leq 2\}$. The random variable $U$ is said to be sub-Gaussian if $\|U\|_{\psi_2}<\infty$ and sub-exponential if $\|U\|_{\psi_1}<\infty$. If $V$ is another random variable that is equal in distribution to $U$, then we write $V \,\scriptstyle\overset{\mathcal{L}}{=}\, \textstyle U$.  In the case when $V$ is a random vector, its correlation matrix is denoted by $\text{cor}(V)$. The standard Gaussian distribution function and density are respectively denoted by $\Phi$ and $\phi$.   If $x$ is a vector, matrix, or tensor with real entries, we use $\|x\|_{\infty}$ to refer to the maximum absolute value of the entries, and $\|x\|_1$ to refer to the sum of the absolute values of the entries. The identity matrix in $\R^{p\times p}$ is denoted by $I_p$. Throughout the paper, the symbol $c$ will denote a positive absolute constant whose value may vary at each occurrence. Different symbols will be used when it is necessary to track constants.
Also, in order to simplify presentation, we define $\Log(t)=\max\{\ln(t),1\}$ for any $t>0$, where $\ln$ is the ordinary natural logarithm.

The class of all hyperrectangles in $\R^p$ is henceforth denoted as $\mathscr{R}$. More precisely, the class $\mathscr{R}$ consists of all cartesian products of intervals $\prod_{j=1}^p \mathcal{I}_j$, where for each $j=1,\dots,p$, the interval $\mathcal{I}_j\subset\R$ is of an arbitrary form (i.e.~it may include zero, one, or two endpoints, and may be bounded or unbounded). 
For any set $A\subset \R^p$ and any positive scalar $t$, an outer $t$-neighborhood of $A$ (in the $\ell_{\infty}$ sense) is defined as 
$A^t=\big\{x\in\R^p\, |\, d(x,A)\leq t\},$
 where the distance from a point to a set is given by $d(x,A)=\inf\{\|x-y\|_{\infty} \, | \, y\in A\}$. In addition, an inner $t$-neighborhood is defined as
 $A^{-t}=\big\{x\in A\, |\, B(x,t)\subset A\},$
 where $B(x,t)=\{y\in\R^p| \|x-y\|_{\infty}\leq t\}$. Furthermore, an associated boundary set of ``width'' $2t$ is defined as $\partial A(t)=A^t\setminus A^{-t}$.\\[-0.2cm]

\noindent\emph{Outline.}
After the main results are presented in Section~\ref{sec:main}, a high-level proof of the central Gaussian approximation result (Theorem~\ref{thm:main}) is given in Section~\ref{sec:mainproof}. Next, some preparatory items are developed in Section~\ref{sec:prep}, which will be used in the more technical arguments given in Section~\ref{sec:delta}. Later on, the other main results (Theorems~\ref{thm:comp}-\ref{thm:bootemp}) are proven in Sections~\ref{sec:proofcomp}-\ref{sec:bootempproof}, and various background results are summarized in Section~\ref{sec:background}.\\[-0.2cm]

\section{Main results}\label{sec:main} The following theorem is the core result of the paper. Subsequently, several other results are obtained as extensions, including a generalization to the case of degenerate covariance matrices (Theorem~\ref{thm:rho0}), a Gaussian comparison result (Theorem~\ref{thm:comp}), and two bootstrap approximation results (Theorems~\ref{thm:boot} and~\ref{thm:bootemp}).

\begin{theorem}[Gaussian approximation]\label{thm:main}
There is an absolute constant $C>0$, such that the following holds for all $n, p\geq 1$ and  $q\in \{1,2\}$: Let $X_1,\dots,X_n\in\R^p$ be centered i.i.d.~random vectors, and suppose that $\nu_q=\max_{1\leq j\leq p}\|X_{1j}/\sqrt{\var(X_{1j})}\|_{\psi_q}$ is finite. In addition, let $\rho$ be the smallest eigenvalue of the correlation matrix of $X_1$, and suppose that $\rho>0$. Lastly, let $Y\in\R^p$ be a centered Gaussian random vector with $\E[YY\ttop]=\E[X_1X_1\ttop]$.
Then, 
\begin{equation}\label{eqn:thm}
\sup_{A\in\mathscr{R}}\Big|\P\big(\ts\frac{1}{\sqrt n}\sum_{i=1}^n X_i\in A\big)-\P(Y\in A)\Big|
  \ \leq \ \displaystyle\frac{ C\nu_q^4\Log^{r_q}(pn)\,\Log(n)}{\rho^{3/2}\, n^{1/2}},
\end{equation}
where $r_1=6$ and $r_2=4$.
\end{theorem}
\noindent \emph{Remarks.}  Although the focal point of this bound is the near $n^{-1/2}$ dependence on $n$,
 it is of interest to assess the dependence on the parameters $\log(p)$, $\nu_q$, and $\rho$ in relation to prior work. For instance, the paper~\citep{Fang:Koike:2020} achieves better dependence on these parameters, with bounds proportional to $\rho^{-1}\log^{3/2}(p)\log(n)n^{-1/2}$ and $(\nu_2/\rho)^{2/3}\log^{4/3}(pn)n^{-1/3}$, corresponding to the respective cases where the data have a log-concave density or sub-Gaussian entries. In Theorem~\ref{thm:main}, the dependence on $\log(p)$ arises mainly from three sources: bounds on the random variable $\|X_1\|_{\infty}^3$ (Lemma~\ref{lem:tail}), bounds on the third-order derivatives of a smoothed indicator function~(Lemma~\ref{lem:deriv}), and Nazarov's Gaussian anti-concentration inequality~(Lemma~\ref{lem:anti}). Meanwhile, with regard to optimal dependence on $\Log(p)$, the paper~\citep{Fang:Koike:2020} provides a lower bound for the left side of~\eqref{eqn:thm}, showing that joint dependence on $(n,p)$ of the form $\log^{3/2}(p)n^{-1/2}$ is generally unimprovable when the data have sub-exponential entries. In addition, the paper~\citep{Das:Lahiri:2020} shows that under certain moment assumptions, the left side of~\eqref{eqn:thm} can approach 0 asymptotically if $\log^2(p)=o(n)$, and that without this condition, convergence to 0 can fail in certain cases.

As for the parameter $\rho$, the assumption of its positivity should be noted, since the paper~\citep{CCKK:2019} establishes the bound $\nu_2^{1/2}\log^{5/4}(pn)n^{-1/4}$ while allowing $\rho=0$. In Theorem~\ref{thm:main}, the positivity assumption is of a technical nature, because it facilitates the implicit smoothing technique. That is, when $\rho>0$, the effect of convolution smoothing from Gaussian partial sums can be more readily quantified. This connection also partly explains the $\rho^{-3/2}$ dependence, because if a univariate function is smoothed by convolution with a $N(0,\rho)$ distribution, then each time the resulting function is differentiated, an extra factor of $\rho^{-1/2}$ is introduced. Hence, in the present context, the bounds on certain third-order remainder terms involve a factor of $\rho^{-3/2}$.
Nevertheless, the next result shows that the case $\rho=0$ can be handled, provided that the correlation matrix of $X_1$ is close to a positive definite matrix in an entrywise sense. Accordingly, in this extended version, the distinct symbol $\varrho$ is used to denote the smallest eigenvalue of the nearby correlation matrix that is positive definite.

As a matter of notation, the quantity $r_q$ will remain as defined in Theorem~\ref{thm:main} throughout the results below.

\begin{theorem}[Gaussian approximation allowing degeneracy]\label{thm:rho0}
There is an absolute constant $c>0$, such that the following holds for all $n,p\geq1$ and $q\in\{1,2\}$:  Let \smash{$X_1,\dots,X_n\in\R^p$} be centered i.i.d.~random vectors with common covariance matrix $\Sigma^X$,
and let $Y$ be a centered Gaussian random vector in $\R^p$ with covariance matrix $\Sigma^Y$. Suppose that $\omega_q=\max_{1\leq j\leq p}\|X_{1j}/\sqrt{\var(Y_{j})}\|_{\psi_q}$ is finite, and that the smallest eigenvalue $\varrho$ of the correlation matrix of Y satisfies $\varrho >0$. Lastly, define $D=\textup{diag}(\Sigma_{11}^Y,\dots,\Sigma_{pp}^Y)$, as well as
\begin{equation*}
%\label{eqn:cltdeltadef}
\Delta=\|D^{-1/2}(\Sigma^X-\Sigma^Y)D^{-1/2}\|_{\infty}.
\end{equation*}
Then,
\begin{equation*}
\sup_{A\in\mathscr{R}}\Big|\P\big(\ts\frac{1}{\sqrt n}\sum_{i=1}^n X_i\in A\big)-\P(Y\in A)\Big|
  \ \leq \ \displaystyle\frac{ c\, \omega_q^4 \,\Log^{r_q}(pn)\,\Log(n)}{\varrho^{3/2}\, n^{1/2}} \ +
  \  \displaystyle \big(\ts\frac{c}{\varrho}\big)\Log(p)\Log(n) \Delta.
\end{equation*}
\end{theorem}

\noindent\emph{Remarks.}   The statement of this result remains true if $\omega_q$ is replaced with $\nu_q$, and  in the context of Theorem~\ref{thm:main} these two parameters are the same. However, the use of $\omega_q$ seems to be somewhat more convenient when applying Theorem~\ref{thm:rho0} to develop results on bootstrap approximation. 

Although Theorem~\ref{thm:rho0} shows that it is possible to achieve near $n^{-1/2}$ rates when $\rho=0$, this is only guaranteed when a choice of $\Sigma^Y$ exists such that the quantities $1/\varrho^{3/2}$ and $n^{1/2}\Delta /\varrho$ are not too large (e.g. polylogarithmic in $n$ and $p$). As an example of a situation where these conditions can be applicable,  our analysis of the bootstrap will use Theorem~\ref{thm:rho0} by letting a (possibly non-invertible) sample covariance matrix play the role of $\Sigma^X$, and letting a population covariance matrix play the role of $\Sigma^Y$. With regard to more general situations, Theorem~\ref{thm:rho0} raises an interesting open question: To what extent can near $n^{-1/2}$ rates of Gaussian approximation be achieved in high dimensions without conditions involving $\varrho$ and $\Delta$? Furthermore, this question can also be asked for bootstrap approximation, since our results in that context have a form that is analogous to Theorem~\ref{thm:rho0}.

\subsection{Gaussian comparison}\label{sec:comp} The next result compares the distributions of two Gaussian random vectors in terms of the entries of their covariance matrices. Whereas techniques related to Stein's method were recently used to establish this result in~\citep[][Theorem 1.1]{Fang:Koike:2020}, it turns out that the result can be obtained in a considerably different way as modification of the proof of Theorem~\ref{thm:main}---which will be shown in Section~\ref{sec:proofcomp}. Later on, this result will serve as a bridge to connect our Gaussian and bootstrap approximation results.

\begin{theorem}[Gaussian comparison]\label{thm:comp}
There is an absolute constant $c>0$, such that the following holds for all $p\geq 1$: Let $Y$ and $Z$ be centered Gaussian random vectors in $\R^p$ having respective covariance matrices $\Sigma^Y$ and $\Sigma^Z$. In addition, let $\varrho$ be the smallest eigenvalue of the correlation matrix of $Y$, and suppose that $\varrho>0$. Lastly, let $D=\textup{diag}(\Sigma_{11}^Y,\dots,\Sigma_{pp}^Y)$, and let
\begin{equation}\label{eqn:gdeltadef}
\Delta^{\!\prime}=\|D^{-1/2}(\Sigma^Z-\Sigma^Y)D^{-1/2}\|_{\infty}.
\end{equation}
 Then,
\begin{equation}\label{eqn:thm2}
\sup_{A\in\mathscr{R}}\Big|\P\big(Z\in A\big)-\P\big(Y\in A\big)\Big|
 \ \leq \ \displaystyle \big(\ts\frac{c}{\varrho}\big)\Log(p)\Log\big(\ts\frac{1}{\Delta^{\!\prime}}\big) \Delta^{\!\prime}.
\end{equation}
\end{theorem}

\noindent\emph{Remarks.} Although the bound depends on the invertibility of $\Sigma^Y$, it is important to note that the bound does \emph{not} depend on the invertibility of $\Sigma^Z$. This is a key property in the context of bootstrap approximation, where $\Sigma^Z$ will represent a sample covariance matrix that is possibly non-invertible. Another comment to make about Theorem~\ref{thm:comp} is its relation to Corollary 5.1 of the paper~\citep{CCKK:2019}. In the notation used here, that result implies a Gaussian comparison bound proportional to $\log(p)\sqrt{\Delta^{\!\prime}}$, which has advantages insofar as it does not depend on $\varrho$, and allows both $\Sigma^{Z}$ and $\Sigma^Y$ to be non-invertible. On the other hand, the current bound has the favorable property that its dependence on $\Delta^{\!\prime}$ is nearly linear.

\subsection{Bootstrap approximation}\label{sec:boot}
Let $X_1^*,\dots,X_n^*$ be sampled with replacement from centered observations $X_1,\dots,X_n$, and let $\bar X=\frac{1}{n}\sum_{i=1}^n X_i$. In this notation, the non-parametric (empirical) bootstrap method of~\cite{Efron:1979} seeks to approximate the distribution of the sum $\frac{1}{\sqrt n}\sum_{i=1}^n X_i$ with the conditional distribution of the sum $\frac{1}{\sqrt n}\sum_{i=1}^n (X_i^*-\bar X)$ given $X_1,\dots,X_n$.  As an alternative to sampling with replacement, the Gaussian multiplier bootstrap method of~\cite{CCK:2013} generates a set of independent random vectors $X_1^{\star},\dots,X_n^{\star}$ from the Gaussian distribution $N(0,\hat\Sigma)$, where we define the sample covariance matrix
\begin{equation}\label{eqn:sigmahatdef}
\hat\Sigma=\frac{1}{n}\displaystyle\sum_{i=1}^n (X_i-\bar X)(X_i-\bar X)\ttop.
\end{equation}
When discussing either method, we will use the notation $\P(\cdot \,| X)$ to refer to probability that is conditional on $X_1,\dots,X_n$.

\begin{theorem}[Multiplier bootstrap approximation]\label{thm:boot}
There is an absolute constant \smash{$c>0$,} such that the following holds for all $n,p\geq 1$ and $q\in\{1,2\}$: Suppose that the conditions of Theorem~\ref{thm:rho0} hold, and let $(\omega_q,\varrho,\Delta)$ be as defined there. In addition, let $X_1^{\star},\dots,X_n^{\star}$ be independent~Gaussian random vectors drawn from $N(0,\hat\Sigma)$.
 Then, the event below holds with probability at least $1-\frac{c}{n}$,
\begin{equation}\label{eqn:thm:boot:mult}
\begin{split}
\sup_{A\in\mathscr{R}}\bigg|\P\Big(\ts\frac{1}{\sqrt n}\sum_{i=1}^n X_i^{\star}\in A\Big|X\!\Big)-\P\Big(\ts\frac{1}{\sqrt n}\sum_{i=1}^n X_i\in A\Big)\bigg|
 & \ \leq \  
 \displaystyle\frac{ c\omega_q^4\Log^{r_q}(pn)\,\Log(n)}{\varrho^{3/2}\, n^{1/2}}\\[0.2cm]
 & \ \ \ \ \ \ \ + \ \displaystyle \big(\ts\frac{c}{\varrho}\big)\Log(p)\Log(n) \Delta.
  \end{split}
\end{equation}
\end{theorem}

\noindent\emph{Remarks.} This result is essentially a direct consequence of Theorem~\ref{thm:rho0} and Theorem~\ref{thm:comp}. There are nevertheless a few details involved, which are given in Section~\ref{sec:bootproof}. 

Regarding the quantities $\varrho$ and $\Delta$, their roles in the bound~\eqref{eqn:thm:boot:mult} may not be intuitive at first sight, since they depend on the correlation matrix of the random vector $Y$ introduced in Theorem~\ref{thm:rho0}. The important point to notice is that the bound~\eqref{eqn:thm:boot:mult} holds for \emph{any choice} of the matrix $\text{cor}(Y)$, provided that its smallest eigenvalue $\varrho$ is positive. So, in other words, the bound involves a theoretical tradeoff between $\varrho$ and $\Delta$ that is governed by the choice of $\text{cor}(Y)$. When $\text{cor}(X_1)$ is positive definite, the simplest choice to consider is $\text{cor}(Y)=\text{cor}(X_1)$, which yields $\varrho=\rho$ and $\Delta=0$. Similarly, these considerations also apply to the next result on the empirical bootstrap.

\begin{theorem}[Empirical bootstrap approximation]\label{thm:bootemp}
There is an absolute constant \smash{$c>0$,} such that the following holds for all $n,p\geq 1$ and $q\in\{1,2\}$: Suppose that the conditions of Theorem~\ref{thm:rho0} hold, and let $(\omega_q,\varrho,\Delta)$ be as defined there. In addition, let $X_1^*,\dots,X_n^*$ be sampled with replacement from $X_1,\dots,X_n$.
 Then, the following event holds with probability at least $1-\frac{c}{n}$,
\begin{equation}\label{eqn:thm:boot}
\begin{split}
\sup_{A\in\mathscr{R}}\bigg|\P\Big(\ts\frac{1}{\sqrt n}\sum_{i=1}^n (X_i^*-\bar X)\in A\Big|X\!\Big)-\P\Big(\ts\frac{1}{\sqrt n}\sum_{i=1}^n X_i\in A\Big)\bigg|
 & \ \leq \  
 \displaystyle\frac{ c\,\omega_q^{4}\,\Log^{s_q}(pn)\,\Log(n)}{\varrho^{3/2}\, n^{1/2}}\\[0.2cm]
 &  \ \ \ \ \ \ \ + \ \big(\ts\frac{c}{\varrho}\big)\Log(p)\Log(n) \Delta,
  \end{split}
\end{equation}
\normalsize
where $s_1=10$ and $s_2=6$.
\end{theorem}

\noindent\emph{Inference on mean vectors}.
From the standpoint of applications, rates of bootstrap approximation for $S_n$ are particularly relevant for doing inference on high-dimensional mean vectors---either via simultaneous confidence intervals  or hypothesis testing.
To illustrate these connections, let $\xi_1,\dots,\xi_n\in\R^p$ be i.i.d.~observations with coordinate-wise means $\mu_j=\E[\xi_{1j}]$ for $j=1,\dots,p$. 
Also define the centered variables $X_{ij}=\xi_{ij}-\mu_j$, as well as the max statistic $M_n=\max_{1\leq j\leq p}\frac{1}{\sqrt n} |S_{nj}|$, again with $S_n=\frac{1}{\sqrt n}\sum_{i=1}^n X_i$. 

The quantiles of $M_n$ have special importance, since they can be used to formulate a set of theoretical confidence intervals for $\mu_1,\dots,\mu_p$. More specifically, if $\bar X_j=\frac{1}{n}\sum_{i=1}^n X_{ij}$, and if $q_{1-\alpha}$ denotes the $1-\alpha$ quantile of $M_n$ for $\alpha\in (0,1)$, then intervals defined by \smash{$\mathcal{I}_j=[\bar X_j \pm \ts\frac{q_{1-\alpha}}{n^{1/2}}]$}
have the property that they simultaneously cover $\mu_1,\dots,\mu_p$ with probability at least $1-\alpha$. Alternatively, in the context of hypothesis testing, the quantile $q_{1-\alpha}$ leads to a theoretical test that rejects the null hypothesis $\mu_1=\cdots=\mu_p=0$ at level $\alpha$. Namely, the test rejects if any of the intervals $\mathcal{I}_1,\dots,\mathcal{I}_p$ exclude 0.

Since the quantile $q_{1-\alpha}$ is unknown, it must be estimated. For this purpose, let $\xi_1^*,\dots,\xi_n^*$ be sampled with replacement from $\xi_1,\dots,\xi_n$, let $\bar\xi=\frac{1}{n}\sum_{i=1}\xi_i$, and define the bootstrapped sum   $S_n^*=n^{-1/2}\sum_{i=1}^n (\xi_i^*-\bar \xi)$, as well as the bootstrapped max statistic $M_n^*=\max_{1\leq j\leq p}|S_{nj}^*|$. Then, an estimate for $q_{1-\alpha}$ can be defined as\label{page:quantile}
$$\hat q_{1-\alpha}= \inf\big\{t\in\R \ \big| \ \P(M_n^*\leq t|X)\geq 1-\alpha\big\}.$$
An estimate based on the Gaussian multiplier bootstrap is defined analogously. 

The coverage quality of $\hat q_{1-\alpha}$ can be analyzed with the previous theoretical results, using the fact that the following relation holds for any hypperectangle of the form $A=[-t,t]^p$ with $t\geq 0$,
\vspace{-0.3cm}
$$\P(M_n^*\leq t|X) \ = \ \P(S_n^*\in A|X).$$
Based on this relation and the previous results, it follows that the estimate $\hat q_{1-\alpha}$ has accurate coverage in the following sense:
If the random vectors $X_1,\dots,X_n$ satisfy the conditions of Theorem~\ref{thm:main}, then the bound
\begin{equation*}
\big|\P(M_n \leq \hat q_{1-\alpha})-(1-\alpha)\big| \ \leq \  \frac{ c\,\nu_q^{4}\,\Log^{s_q}(pn)\,\Log(n)}{\rho^{3/2}\, n^{1/2}}.
\end{equation*}
holds for some absolute constant $c>0$. (See Lemma~\ref{lem:quantiles} for details.)
For a more extensive overview of the links between bootstrap approximation results and high-dimensional inference, we refer to the survey~\citep{Belloni:2018}, which discusses the ``many approximate means'' framework, and includes an abundance of further references.

\section{The main steps in the proof of Theorem~\ref{thm:main}}\label{sec:mainproof}

Before starting the proof, we first comment on some details and basic reductions that will be used throughout the remainder of the paper, often without explicit mention. \\[-0.2cm]

\noindent\emph{Preliminary comments.}\label{page:prelim} First, note the inequalities $\rho,\varrho \leq 1$ and $\nu_q\geq 1$, which follow from the definitions of $\rho$, $\varrho$, and $\nu_q$. Next, for the theorems involving the quantity $\Delta$, we may assume without loss of generality that \smash{$\Delta\leq 1/2$,} for otherwise such results are true when the absolute constant $c$ is taken to be at least $2$. Under this reduction, we have $1/2 \ \leq \var(X_{1j})/\var(Y_j)\leq 3/2$  for all $j=1,\dots,p$, which implies that $\omega_q$ and $\nu_q$ are interchangeable up to absolute constants: $\sqrt{1/2}\,\nu_q\leq \omega_q \leq \sqrt{3/2}\,\nu_q$. From this, it is also useful to note that $\omega_q\geq \sqrt{1/2}$. Lastly, we may assume without loss of generality, that the covariance matrix $\E[YY\ttop]$ has all ones along the diagonal, because the distributional metric we use is invariant to diagonal rescaling. Hence, we may write $\Delta=\|\Sigma^X-\Sigma^Y\|_{\infty}$ and $\Delta^{\!\prime}=\|\Sigma^Z-\Sigma^Y\|_{\infty}$.\\[-0.2cm]

\noindent\emph{Lindeberg interpolation.} To lay the groundwork for our use of Lindeberg interpolation, let \smash{$Y_1,\dots,Y_n\in\R^p$} be i.i.d.~copies of $Y$ that are independent of $X_1,\dots,X_n$. When dealing with various partial sums, we use the following notation for $1\leq k\leq k'\leq n$,
\begin{equation*}
\begin{split}
S_{k:k'}(X)&=n^{-1/2}(X_k+\dots+X_{k'})\\[0.2cm]
S_{k:k'}(Y)&=n^{-1/2}(Y_k+\cdots+Y_{k'}).
\end{split}
 \end{equation*}
In addition, it will be convenient to denote the distance between $S_{1:k}(X)$ and $S_{1:k}(Y)$ as
 $$\mathsf{D}_k =\sup_{A\in\mathscr{R}}\Big|\P(S_{1:k}(X)\in A)-\P(S_{1:k}(Y)\in A)\Big|.$$
The interpolation for bounding $\mathsf{D}_n$ is based on the following quantities, which are defined for each $A\in\mathscr{R}$ according to
 \begin{align*}
 \delta_k^X(A) & \ =\ \P\Big(S_{1:k-1}(X)+\ts\frac{1}{\sqrt n} X_k+S_{k+1:n}(Y)\in A\Big) \ -\ \P\Big(S_{1:k-1}(X)+S_{k+1:n}(Y)\in A\Big),\\[0.2cm]
 \delta_k^Y(A) &\ = \  \P\Big(S_{1:k-1}(X)+\ts\frac{1}{\sqrt n} Y_k+S_{k+1:n}(Y)\in A\Big) \ \ - \ \,\P\Big(S_{1:k-1}(X)+S_{k+1:n}(Y)\in A\Big).
 \end{align*}
This notation yields the interpolation
\begin{equation*}
\P(S_{1:n}(X)\in A)-\P(S_{1:n}(Y)\in A) \ = \ \sum_{k=1}^n  \{\delta_k^X(A)-\delta_k^Y(A)\}.
\end{equation*}
Next, define the supremum of the $k$th absolute difference
 \begin{equation}\label{eqn:deltak}
  \delta_k \ = \ \sup_{A\in\mathscr{R}}|\delta_k^X(A)-\delta_k^Y(A)|,
  \end{equation}
  which leads to the bound
  $$\mathsf{D}_n \ \leq \  \delta_1+\cdots+\delta_n.$$
However, rather than working directly with the entire sum $\delta_1+\cdots+\delta_n$, we will begin with a lemma  that reduces the problem to bounding $\delta_1+\cdots+\delta_{n-m}$ for an integer $m$ that will be chosen later.

\begin{lemma}\label{lem:start}
There is an absolute constant $c_1>0$ such that the following holds for all $n,p\geq 1$, $m\in\{1,\dots,n\}$, and $q\in\{1,2\}$:
If the conditions of Theorem~\ref{thm:main} hold, then
\begin{equation}\label{eqn:lem:start}
\mathsf{D}_n \ \leq \ c_1\nu_q \sqrt{\ts\frac{m}{n}}(\Log(pn))^{\frac{2+q}{2q}} \ + \ 3(\delta_1+\cdots+\delta_{n-m}).
\end{equation}
\end{lemma}
 \noindent \emph{Remark.} This statement is also true under the conditions of Theorem~\ref{thm:rho0}, provided that $\nu_q$ is replaced with $\omega_q+1$. Here, it is worthwhile to proceed straight to the proof of this lemma, since the argument is fairly short.

\proof As a temporary shorthand, let $\zeta=S_{1:n}(X)$ and $\xi=S_{1:n-m}(X)+S_{n-m+1:n}(Y)$. Observe that the distance between $\xi$ and $S_{1:n}(Y)$ is at most $\delta_1+\cdots+\delta_{n-m}$, which gives
\begin{equation*}
\mathsf{D}_n \ \leq \ \sup_{A\in\mathscr{R}}\Big|\P\big(\zeta\in A\big) -\P\big(\xi\in A\big)\Big| \ + \ (\delta_1+\cdots+\delta_{n-m}).
\end{equation*}
To control the first term on the right side, Lemma~\ref{lem:coupling} gives the following basic inequality, which holds for any $A\in\mathscr{R}$ and any $t>0$,
\begin{equation*}
|\P(\zeta\in A) -\P(\xi\in A)| \ \leq \ \P\big(\xi \in \partial A(t)\big) \ + \ \P(\|\zeta-\xi\|_{\infty}\geq t),
\end{equation*}
where the boundary set $\partial A(t)$ is defined in the notation paragraph of Section~\ref{sec:intro}. 
Next, we control the probability that $\xi$ hits $\partial A(t)$ by replacing $\xi$ with $S_{1:n}(Y)$ at the price of some error terms. To do this, again note that the distance between $\xi$ and $S_{1:n}(Y)$ is at most \smash{$\delta_1+\cdots+\delta_{n-m}$,} and so
\begin{equation*}
\begin{split}
\P\big(\xi \in A^t\big) 
& \ \leq \ \P\Big( S_{1:n}(Y)  \in A^t\Big) \, + \, (\delta_1+\cdots+\delta_{n-m}).
\end{split}
\end{equation*}
Similarly
\begin{equation*}
\begin{split}
\P\big(\xi \in A^{-t}\big) 
& \ \geq \ \P\Big( S_{1:n}(Y)  \in A^{-t}\Big) - (\delta_1+\cdots+\delta_{n-m}),
\end{split}
\end{equation*}
and then combining gives
\begin{equation*}
 \P\big(\xi \in \partial A(t)\big) \ \leq \ \P\Big(S_{1:n}(Y)  \in \partial A(t)\Big) \ + \ 2(\delta_1+\cdots+\delta_{n-m}).
\end{equation*}
In turn, Nazarov's Gaussian anti-concentration inequality (Lemma~\ref{lem:anti}) gives
\begin{equation*}
\P\Big(S_{1:n}(Y)  \in \partial A(t)\Big)  \ \leq 
\  ct\sqrt{\Log(p)},
\end{equation*}
where we have made use of the reduction that $\E[YY\ttop]$ has all ones along the diagonal.
Finally, the tail bound in part (iv) of Lemma~\ref{lem:tail} implies that if 
$t=c\nu_q\sqrt{\ts\frac{m}{n}}\Log^{1/q}(pn)$ for a sufficiently large absolute constant $c>0$, then the probability $\P(\|\zeta-\xi\|_{\infty}\geq t)$ is at most $\frac{c}{n}$.\qed

 ~\\[-0.2cm]

\noindent\emph{Comments on induction.}
 In order to prove Theorem~\ref{thm:main}, we will use a form of strong induction. Specifically, for a given absolute constant $C>0$, and given integers $n$ and $p$, the associated induction hypothesis is that the inequality below (denoted  $\textup{H}_k(C)$) holds simultaneously for all $k=1,\dots,n-1$, 
 \begin{equation}\label{eqn:hyp}
 \textup{H}_k(C): \quad \quad \mathsf{D}_k \ \leq \ \frac{C \,\nu_q^{4}\,\Log^{r_q}(pk)\,\Log(k)}{\rho^{3/2} \, k^{1/2}}.
 \end{equation}
 Although it is common in high-dimensional statistics to think of $n$ and $p$ as growing together, it is worth clarifying that the inductive approach here is based on showing that, for any fixed $p$, the entire sequence $\textup{H}_1(C),\textup{H}_2(C),\dots$ holds.  Hence, because $p$ is arbitrary, it will follow that the statement of the main result holds for all pairs $(n,p)$.

\subsection{Proof of Theorem~\ref{thm:main}} Observe that if $C\geq \sqrt 2$, then it is clear that $\textup{H}_1(C)$ and $\textup{H}_2(C)$ hold. To carry out the induction, fix any $n\geq 3$, and suppose that $\textup{H}_1(C),\dots,\textup{H}_{n-1}(C)$ hold for some absolute constant $C\geq \sqrt 2$. Our goal is now to show that $\textup{H}_n(C)$ holds with the same value of $C$. The main tool for this purpose is the proposition below, whose proof is deferred to Section~\ref{sec:delta}.

\begin{proposition}\label{prop:sum} There is a positive absolute constant $c_2$ such that the following holds for all $n\geq 3$, $p\geq 1$, $ m\in [1, n/3]$, and $q\in\{1,2\}$: If the conditions of Theorem~\ref{thm:main} hold, and if $\textup{H}_1(C),\dots,\textup{H}_{n-1}(C)$ hold for some absolute constant $C\geq \sqrt 2$, then
\begin{equation}\label{eqn:prop:sum}
\delta_1+\cdots+\delta_{n-m} \ \leq
 \ \frac{\big(\ts\frac{c_2\nu_q^{3}}{\rho}\big)(\Log(pn))^{\frac{5q+6}{2q}}\Log(n)}{n^{1/2}} 
\ + \  \frac{ \big(\ts\frac{c_2 C\nu_q^7}{\rho^3}\big)(\Log(pn))^{r_q+\frac{3q+6}{2q}}\Log(n)}{n^{1/2}m^{1/2}}.
\end{equation}
\end{proposition}

 \noindent \emph{Remark.} At a high level, Lemma~\ref{lem:start} and Proposition~\ref{prop:sum} reduce the proof of Theorem~\ref{thm:main} to exhibiting suitable values of $m$ and $C$, as presented below. Since the details of finding such values are of minor conceptual importance, the reader may wish to skip ahead to Section~\ref{sec:prep}, which covers the more significant aspects of the proof related to implicit smoothing and moment matching.\\[-0.2cm]

To proceed, let $c_1$ and $c_2$ be the absolute constants in the statements of  Lemma~\ref{lem:start} and Proposition~\ref{prop:sum}. We may assume without loss of generality that $c_2=c_1$ and $c_1\geq 1$, because these results remain true if $c_1$ and $c_2$ are both replaced by $\max\{c_1,c_2,1\}$. Next, define the quantities $\alpha$, $\beta$, and $\gamma$ according to
\begin{align*}
\alpha  \ & = \  \frac{c_1\nu_q m^{1/2}(\Log(pn))^{\frac{2+q}{2q}}}{n^{1/2}}\\[0.3cm]
\beta \ & = \   \ \frac{\big(\ts\frac{c_1\nu_q^{3}}{\rho}\big)(\Log(pn))^{\frac{5q+6}{2q}}\Log(n)}{n^{1/2}} 
  \ + \  
\  \frac{ \big(\ts\frac{c_1 C\nu_q^7}{\rho^3}\big)(\Log(pn))^{r_q+\frac{3q+6}{2q}}\Log(n)}{n^{1/2}m^{1/2}}\\[0.3cm]
\gamma \ & = \  \frac{ \big(\frac{C\nu_q^4}{\rho^{3/2}}\big)\Log^{r_q}(pn)\,\Log(n)}{ n^{1/2}}.
\end{align*}
In terms of this notation, Lemma~\ref{lem:start} and Proposition~\ref{prop:sum} give the bound
%\begin{equation}
$\mathsf{D}_n \ \leq \ \alpha \ + \ 3\beta$.
%\end{equation}
Therefore, in order to show $\textup{H}_n(C)$, it is enough to show that there exist choices of $m$ and $C$ for which
\begin{equation}\label{eqn:abc}
\alpha \ + \ 3\beta \ \leq \ \gamma.
\end{equation}
(It is not immediately obvious that such choices exist, because both sides of~\eqref{eqn:abc} depend on $C$, and also, because $m$ must simultaneously satisfy the constraint $m\leq n/3$.) 

In the remainder of the proof, we will construct feasible choices of $C$ and $m$ explicitly in terms of $c_1$.
For this purpose, let $\kappa\geq 1$ be a value that will be tuned later, and consider a choice of $m$ whose square root given by
\begin{equation}\label{eqn:mkappachoice}
m^{1/2} \ = \ \Big\lceil \kappa \big(\ts\frac{\nu_q^{3}}{\rho^{3/2}}\big)\big(\Log(pn)\big)^{\frac{qr_q+2+q}{2q}}\Big\rceil.
\end{equation}
When $m$ is chosen this way, the quantities $\alpha$ and $3\beta$ satisfy
\begin{align}
\alpha  \ & \leq \  \frac{ 2\kappa c_1\big(\ts\frac{\nu_q^{4}}{\rho^{3/2}}\big) (\Log(pn))^{\frac{qr_q+2q+4}{2q}}}{n^{1/2}}\label{eqn:alphatemp} \ \ \ \ \ \ \ \  \\[0.3cm]
3\beta \ & \leq \   \frac{3c_1\big(\ts\frac{\nu_q^{3}}{\rho}\big)(\Log(pn))^{\frac{5q+6}{2q}}\Log(n)}{n^{1/2}}  \ + \  
\frac{\big(\frac{3c_1C}{\kappa}\big) \big(\ts\frac{\nu_q^{4}}{\rho^{3/2}}\big)\big(\Log(pn)\big)^{\frac{qr_q+2q+4}{2q}}\Log(n)}{n^{1/2}}.
\end{align}
At this stage it is helpful to note some basic details. First, we have $\frac{\nu_q^{3}}{\rho}\leq \frac{\nu_q^{4}}{\rho^{3/2}}$,
since $\nu_q\geq 1$ and $\rho\leq 1$. Second, the definitions $r_1=6$ and $r_2=4$ imply that $\frac{qr_q+2q+4}{2q}=r_q$ and $\frac{5q+6}{2q}\leq r_q$ for each $q\in\{1,2\}$. With these items in mind, the last several steps imply
\begin{equation*}
 \alpha \ + \ 3\beta \ \leq \ \Big(2\kappa c_1 +3c_1+\big(\ts\frac{3c_1C}{\kappa}\big)\Big)\cdot \displaystyle\frac{\gamma}{C}.
\end{equation*}
Thus, in order to show $\alpha+3\beta\leq \gamma$, it suffices to select $\kappa$ and $C$ in terms of $c_1$ so that
\begin{equation*}
2\kappa c_1 +3c_1+\big(\ts\frac{3c_1C}{\kappa}\big) \ \leq \ C.
\end{equation*}
Likewise, if we put 
\begin{equation}\label{eqn:kappachoice}
 \kappa=\sqrt{\ts\frac{3}{2}C},
\end{equation}
then $C$ should be chosen to satisfy
\begin{equation*}
(2\sqrt{6}c_1)\sqrt{C}+ 3c_1 \leq C.
\end{equation*}
This is a quadratic inequality in $\sqrt{C}$, which holds when
\begin{equation}\label{eqn:c1lower}
C \ \geq \  \Big(\sqrt 6 c_1+ \sqrt{ 6 c_1^2+3c_1}\Big)^2.
\end{equation}
In particular, this is compatible with the condition $C\geq \sqrt 2$ mentioned earlier, since $c_1\geq 1$.
Moreover, since the right side of~\eqref{eqn:c1lower} is purely a function of $c_1$, the only remaining consideration is to make sure that~\eqref{eqn:c1lower} allows for a feasible choice of $m\leq n/3$. (Note that $m$ is now determined by $C$ through~\eqref{eqn:mkappachoice} and~\eqref{eqn:kappachoice}.) To do this, we may assume without loss of generality that the inequality
\begin{equation*}
n^{1/2} \ \geq \ \big(\ts\frac{C\nu_q^{4}}{\rho^{3/2}}\big)\Log^{r_q}(pn)\Log(n)
\end{equation*}
holds, for otherwise $\textup{H}_n(C)$ is true. Comparing this inequality with~\eqref{eqn:mkappachoice} shows that the condition $m\leq n/3$ holds, for instance, when $2\sqrt{(3/2)C} \ \leq \ C/\sqrt 3$, i.e.~when $C\geq 18$. But at the same time, the right side of~\eqref{eqn:c1lower} is already greater than 18, and so it suffices to take $C$ equal to the right side of~\eqref{eqn:c1lower}.\qed

 \section{Preparatory items}\label{sec:prep}

Below, we develop the notation and key objects that will be needed to prove Proposition~\ref{prop:sum} in Section~\ref{sec:delta}. 
 \subsection{Implicit smoothing}\label{sec:smooth}
 The main idea in this subsection is to represent the quantities $\delta_k^X(A)$ and $ \delta_k^Y(A)$ in terms of a certain implicit Gaussian smoothing function. We use the word ``implicit'', because the smoothing function is automatically built into the Lindeberg interpolation through the Gaussian partial sums.
 
 To proceed, let $\zeta \sim N(0,I_p)$ be a standard Gaussian random vector in $\R^p$, and for any fixed \smash{$s\in\R^p$}, $A\in\mathscr{R}$, and $\e>0$, define 
  \begin{equation}\label{eqn:phiprod}
  \varphi_{\e}(s,A)=\P(s+\e \zeta\in A).
  \end{equation}
For the sake of understanding this function, it may be helpful to note that in the particular case when $A$ has the form $A=\prod_{j=1}^p[a_j,b_j]$, we have the explicit formula
    \begin{equation}\label{eqn:phiprod2}
  \varphi_{\e}(s,A)=\prod_{j=1}^p \Big(\Phi\big(\ts\frac{b_j-s_j}{\e}\big)-\Phi\big(\ts\frac{a_j-s_j}{\e}\big)\Big).
  \end{equation}
When $A$ is held fixed, the function $\varphi_{\e}(\cdot,A)$ is a smoothed version of the indicator \smash{$s\mapsto 1\{s\in A\}$,} with $\e$ playing the role of a smoothing parameter.
 Next, for each $k=1,\dots,n-1$, define 
 \begin{equation}\label{eqn:ekdef}
 \e_k= \ts\sqrt{\frac{n-k}{n}}\sqrt \rho.
 \end{equation}
 The parameter $\e_k$ is used in order to simplify the following (distributional) decomposition of the Gaussian random vector $S_{k+1:n}(Y)$,
 \begin{equation*}
 S_{k+1:n}(Y) \ \ \, \scriptstyle\overset{\mathcal{L}}{\displaystyle=} \, \ \ \displaystyle \e_k V_{k+1} \ + \ \ts\sqrt{\frac{n-k}{n}}W_{k+1},
 \end{equation*}
 where $V_{k+1}\sim N(0, I_p)$ and $W_{k+1}\sim N(0,R-\rho I_p)$ are independent, and $R$ is the correlation matrix of $Y$ with smallest eigenvalue $\rho$.  Here, we continue to work under the reduction that $\E[Y_1Y_1\ttop]=R$. Also, note that the vectors $V_{k+1}$ and $W_{k+1}$ may be taken to be independent of $X_1,\dots,X_n$. Consequently, if we let $\hat A_{k+1}\in\mathscr{R}$ denote the randomly shifted version of $A$ defined by
 \begin{equation}\label{eqn:rhatdef}
 \hat A_{k+1} = \Big\{x-\ts\sqrt{\ts\frac{n-k}{n}}W_{k+1}\Big| x\in A\Big\},
 \end{equation}
 then we can connect $\varphi_{\e_k}$ to the partial sums in the Lindeberg interpolation through the following exact relation
\begin{equation*}
\begin{split}
\P\Big(S_{1:k}(X)+S_{k+1:n}(Y)\in A\Big)  %
& \  = \ \E\Big[\varphi_{\e_k}(S_{1:k}(X)\,,\,\hat A_{k+1})\Big].
\end{split}
\end{equation*}
In turn, this relation allows us to express $\delta_k^X(A)$ in terms of $\varphi_{\e_k}$ for $k=1,\dots,n-1$,
\begin{equation}\label{eqn:deltaksmooth}
\delta_k^X(A)=\E\Big[\varphi_{\e_k}\Big(S_{1:k-1}(X)+\ts\frac{1}{\sqrt n} X_k,\hat A_{k+1}\Big)-\varphi_{\e_k}\Big(S_{1:k-1}(X),\hat A_{k+1}\Big)\Big].
\end{equation}
The formula~\eqref{eqn:deltaksmooth} is the key item to take away from the current subsection. The corresponding expression for $\delta_k^Y(A)$ is nearly identical, with the only change being that the single occurrence of $ X_k$ in~\eqref{eqn:deltaksmooth} is replaced with $Y_k.$

 \subsection{Moment matching}\label{sec:moment}
 By expanding the function $\varphi_{\e_k}(\cdot,\hat A_{k+1})$ to second order at the point $S_{1:k-1}(X)$, 
 we have the moment-matching formulas
\begin{align}
\delta_k^X(A) & \ = \ \E[L_k^X(A)]  \ + \  \E[Q_k^X(A)] \ + \ \E[R_k^X(A)]\label{eqn:deltaxexpand}\\[0.3cm]
\delta_k^Y(A) & \ = \ \ \E[L_k^Y(A)] \ + \ \E[Q_k^Y(A)] \ + \ \E[R_k^Y(A)],\label{eqn:deltayexpand}
\end{align}
where the terms corresponding to $\delta_k^X(A)$ are defined as follows. Specifically, if all derivatives are understood as being with respect to the first argument of $\varphi_{\e_k}$, then
\begin{align}
L_k^X(A) & \ = \ \Big\langle \nabla \varphi_{\e_k}(S_{1:k-1}(X),\hat A_{k+1}), \ n^{-1/2}X_k\Big\rangle\label{eqn:linear}\\[0.2cm]
Q_k^X(A) & \ = \ \ts\frac{1}{2}\Big\langle \nabla^2 \varphi_{\e_k}(S_{1:k-1}(X),\hat A_{k+1}), \ n^{-1} X_kX_k\ttop \Big\rangle\label{eqn:quad}\\[0.3cm]
R_k^X(A) & \ = \ \ts\frac{(1-\tau)^2\!\!}{2}\,\Big\langle \nabla^3 \varphi_{\e_k}\Big(S_{1:k-1}(X)+ \ts\frac{\tau}{\sqrt n}  X_k,\hat A_{k+1}\Big), \ n^{-3/2}X_k^{\otimes 3}\Big\rangle,\label{eqn:RXdef}
\end{align}
with $\tau$ being a Uniform[0,1] random variable that is independent of all other random variables. The notation $\nabla^3\varphi_{\e_k}(s,\hat A_{k+1})$ refers to the tensor in $\R^{p\times p\times p}$ whose entries are comprised by all possible three-fold partial derivatives of $\varphi_{\e_k}(\cdot,\hat A_{k+1})$ at the point $s$. Also, we use $\langle \cdot,\cdot\rangle$ to denote the entrywise inner product on vectors, matrices, and tensors.
Lastly, the terms $L_k^Y(A)$, $Q_k^Y(A)$ and $R_k^Y(A)$ associated with $\delta_k^Y(A)$ in~\eqref{eqn:deltayexpand} only differ from those above insofar as each appearance of $X_k$ on the right sides of~\eqref{eqn:linear},~\eqref{eqn:quad}, and~\eqref{eqn:RXdef} is replaced by $Y_k$.

The classical idea of the Lindeberg interpolation is that if~\eqref{eqn:deltayexpand} is subtracted from~\eqref{eqn:deltaxexpand}, then the first and second order terms cancel. This is because $X_k$ and $Y_k$ have matching mean vectors and covariance matrices, and are independent of $S_{1:k-1}(X)$ and $\hat A_{k+1}$.
Consequently, we have the relation
\begin{equation}\label{eqn:deltadiffmatch}
\delta_k^X(A)-\delta_k^Y(A) \ = \ \E[R_k^X(A)]-\E[R_k^Y(A)].
\end{equation}

Hence, in order to control the supremum $\delta_k=\sup_{A\in\mathscr{R}}|\delta_k^X(A)-\delta_k^Y(A)|$ in~\eqref{eqn:deltak}, it remains to bound the expected remainders, and this is handled in the next section. An important technical tool for this task is a result that provides bounds on the partial derivatives of $\varphi_{\e_k}(\cdot,A)$ holding uniformly with respect to $A\in\mathscr{R}$.  The earliest variant of this result we are aware of is~\citep[][Theorem 3]{Bentkus:1990}. Other variants were later developed in~\cite[][eqn.~2.10]{Anderson:1998}, and~\cite[][Lemmas 2.2 and 2.3]{Fang:Koike:2020}. The following version is a slight reformulation of a special case of the last of these results. As a matter of notation, recall that if $x$ is a real tensor, we use $\|x\|_1$ to denote the sum of the absolute values of its entries.
\begin{lemma}\label{lem:deriv}
Fix an integer $r\geq 1$, and let $\e>0$. Then, there is a constant $c_r>0$ depending only on $r$ such that
\begin{equation}
 \sup_{(s,A)\in \mathbb{R}^{^p} \times\mathscr{R}}\|\nabla^r \varphi_{\e}(s,A)\|_1 \ \leq \ \frac{c_r\,\Log^{r/2}(p)}{\e^r}.
\end{equation}
\end{lemma}
\noindent In relation to the early result~\citep[][Theorem 3]{Bentkus:1990}, the main advantage of this version is that it applies to the entire class of hyperrectangles $\mathscr{R}$, rather than to certain sub-classes. See also~\citep [Theorem 6.5]{Fooling:2019} for a more recent formulation of Bentkus' result.

\section{Bounds for $\delta_k$, and the proof of Proposition~\ref{prop:sum}}\label{sec:delta} 
The next lemma handles $\delta_k$ for $k=2,\dots,n-1$. This lemma is of special significance to the overall structure of the proof of Theorem~\ref{thm:main}, because it sets up the opportunity to apply the induction hypothesis to $\mathsf{D}_{k-1}$. Note that because the upper bound in the lemma involves the quantity $1/(k-1)$, the case of $k=1$  will be handled separately in Lemma~\ref{lem:k1} later on. It will not be necessary to handle $\delta_n$, due to Lemma~\ref{lem:start}.  At the end of the section, the proof of Proposition~\ref{prop:sum} will be given.

\begin{lemma}\label{lem:k2}
There is an absolute constant $c>0$ such that the following holds for all $n\geq 3$, $p\geq 1$, and $k\in\{2,\dots,n-1\}$ and $q\in\{1,2\}$:
If the conditions of Theorem~\ref{thm:main} hold, then
\begin{equation}\label{eqn:deltakbound2}
\delta_k \ \leq  \ \ \ts\frac{c\nu_q^{3} (\Log(p))^{\frac{3q+6}{2q}}}{\e_k^3\, n^{3/2}}\Big(\e_k\Log(pn)\sqrt{\ts\frac{n}{k-1}}\, +\, \mathsf{D}_{k-1} \, + \, \ts\frac{1}{pn}\Big).
\end{equation}
\end{lemma}
\proof From the previous section, we have the following bound on $\delta_k$,
\begin{equation}\label{eqn:deltaksplit}
 \delta_k \ \leq \ \sup_{A\in\mathscr{R}}\E[|R_k^X(A)|] \ + \ \sup_{A\in\mathscr{R}}\E[|R_k^Y(A)|].
\end{equation}
The current proof will only establish a bound on $\sup_{A\in\mathscr{R}}\E[|R_k^X(A)|]$, since the same argument can be applied to $\sup_{A\in\mathscr{R}}\E[|R_k^Y(A)|]$. To begin, let $\tilde A_{k+1}$ be the random hyperrectangle obtained by subtracting $\ts\frac{\tau}{\sqrt n}X_k$ from all points in $\hat A_{k+1}$,
and for any fixed $\ve>0$, define the event
$$E_k(\ve) \ = \ \Big\{ S_{1:k-1}(X) \in \partial \tilde A_{k+1}(\ve)\Big\}.$$
Below, we will separately analyze $R_k^X(A)$ on the event $E_k(\ve)$ and its complement $E_k(\ve)^c$.
~\\[-0.6cm]

\noindent\emph{Handling the remainder on $E_k(\ve)$.}
 By applying H\"older's inequality to the definition of $R_k^X(A)$ in~\eqref{eqn:RXdef}, we have
\begin{equation}\label{eqn:remonak}
|R_k^X(A)| 1\{E_k(\ve)\} \ \leq \ \frac{1}{n^{3/2}}\cdot \Big(\sup_{(s,A)\in\R^{^p}\times \mathscr{R}}\|\nabla^3 \varphi_{\e_k}(s,A)\|_1\Big)\cdot \|X_k\|_{\infty}^3\cdot 1\{E_k(\ve)\},
\end{equation}
Crucially, the second factor on the right can be bounded by Lemma~\ref{lem:deriv}, which gives 
\begin{equation}\label{eqn:bentkus}
\sup_{(s,A)\in\R^{^p}\times \mathscr{R}}\|\nabla^3 \varphi_{\e_k}(s,A)\|_1 \ \leq \ \ts\frac{c\Log^{3/2}(p)}{\e_k^3}.
\end{equation}
Thus, it remains to control the expectation $\E[\|X_k\|_{\infty}^3 1\{E_k(\ve)\}]$. Noting that $S_{1:k-1}(X)$ is independent of $\tilde A_{k+1}$ and $X_{k}$, we have
\begin{equation*}
\begin{split}
\E[\|X_k\|_{\infty}^31\{E_k(\ve)\}] & \ = \ \E\Big[\|X_k\|_{\infty}^3\,\P\Big( S_{1:k-1}(X)\in \partial  \tilde A_{k+1}(\ve)\,\Big|\, \tilde A_{k+1}, X_k\Big)\Big]\\[0.3cm]
& \ \leq \ \E\Big[\|X_k\|_{\infty}^3\,\Big(\P\Big( S_{1:k-1}(Y)\in \partial \tilde A_{k+1}(\ve)\,\Big|\,  \tilde A_{k+1}, X_k\Big) \ + \ 2\mathsf{D}_{k-1}\Big)\Big]\\[0.3cm]
& \ \leq \ \ \E[\|X_k\|_{\infty}^3]\Big(c\ve\sqrt{\ts\frac{n}{k-1}} \sqrt{\Log(p)} \ + \ 2\mathsf{D}_{k-1}\Big),\label{eqn:particprob}
\end{split}
\end{equation*}
where we note that $S_{1:k-1}(X)$ has been replaced with $S_{1:k-1}(Y)$ at the price of $2\mathsf{D}_{k-1}$, and Nazarov's Gaussian anti-concentration inequality (Lemma~\ref{lem:anti}) has been used in the last step.
  Combining the last several steps with the bound $\E[\|X_k\|_{\infty}^3]\leq c\nu_q^3 \Log^{3/q}(p)$ from Lemma~\ref{lem:tail} yields
  \begin{equation}\label{eqn:complicated}
   \E[|R_k^X(A)| 1\{E_k(\ve)\}] \ \leq  \ \ts\frac{c\nu_q^{3}(\Log(p))^{\frac{3q+6}{2q}}}{\e_k^3 n^{3/2}}\Big(\ve\sqrt{\ts\frac{n}{k-1}}\sqrt{\Log(p)}+\mathsf{D}_{k-1}\Big),
  \end{equation}
  which holds uniformly with respect to $A\in\mathscr{R}$.\\[-0.2cm]
 
 \noindent\emph{Handling the remainder on $E_k^c(\ve)$.} For this part, the idea is that for any $A\in\mathscr{R}$, the quantity $\|\nabla^3 \varphi_{\e_k}(s, A)\|_1$ is essentially negligible when $s\not\in\partial A(\ve)$ and $\ve$ is chosen to be sufficiently large. To this end, define the deterministic quantity
$$b_k(\ve)=\sup\Big\{\|\nabla^3\varphi_{\e_k}(s,A)\|_1 \, \Big| \ A\in\mathscr{R}\, \text{ and } \, s\in (\R^p\setminus \partial A(\ve))\Big\},$$
where the supremum involves both $s$ and $A$.
Thus, H\"older's inequality gives
$$\E[|R_k^X(A)|1\{E_k^c(\ve)\}] \ \leq \ \ts\frac{1}{n^{3/2}}\, b_k(\ve) \, \E[\|X_k\|_{\infty}^3].$$
It is shown in Lemma~\ref{lem:bkbound} that if $\ve$ is chosen as $\ve=4\e_k\sqrt{\Log(pn)}$, then 
 \begin{equation}\label{eqn:bkbound}
 b_k(\ve)\leq \ts\frac{c}{\e_k^3pn}.
 \end{equation}
 Combining this with the fact that $\E[\|X_k\|_{\infty}^3]\leq c \nu_q^3\Log^{3/q}(p)$ leads to the stated result.
\qed

 \begin{lemma}\label{lem:k1}
 There is an absolute constant $c >0$, such that the following holds for all $n\geq 2$, $p\geq 1$ and $q\in\{1,2\}$: If the conditions of Theorem~\ref{thm:main} hold, then 
 \begin{equation*}
  \delta_1 \ \leq  \ \ts\frac{c \nu_q^{3}\, (\Log(p))^{\frac{3q+6}{2q}}}{\rho^{3/2}n^{3/2}}.
  \end{equation*}
 \end{lemma}
 \proof As in the proof of the previous lemma, it suffices to bound $\sup_{A\in\mathscr{R}}\E[|R_1^X(A)|]$. Using the same steps as in~\eqref{eqn:remonak} and~\eqref{eqn:bentkus}, but ignoring the role of the indicator $1\{E_k(\ve)\}$, we have
$$\sup_{A\in\mathscr{R}}\E[|R_1^X(A)|] \ \leq \ \ts\frac{c\E[\|X_1\|_{\infty}^3]\Log^{3/2}(p)}{\e_1^3 n^{3/2}}.$$
Applying the previously used bound on $\E[\|X_1\|_{\infty}^3]$ from Lemma~\ref{lem:tail} completes the proof.\qed\\

\noindent\textbf{Proof of Proposition~\ref{prop:sum}.} \ 
By Lemma~\ref{lem:k1}, the quantity $\delta_1$ is negligible in comparison to the right side of~\eqref{eqn:prop:sum}, and so it is enough to focus on $\delta_2+\cdots+\delta_{n-m}$. By Lemma~\ref{lem:k2}, we have that for $k=2,\dots,n-m$,
\begin{equation*}
  \delta_k  \ \leq \ \ts\frac{c\nu_q^{3}(\Log(pn))^{\frac{5q+6}{2q}}}{\e_k^2 n\sqrt{k-1}} 
  \ + \ \frac{c\nu_q^{3} (\Log(pn))^{\frac{3q+6}{2q}}(\mathsf{D}_{k-1}+\frac{1}{pn})}{\e_k^3n^{3/2}}.\\[0.3cm]
\end{equation*}
  Since we assume that $\textup{H}_1(C),\dots,\textup{H}_{n-1}(C)$ hold, we may derive a bound on $\delta_k$ for each $k=2,\dots,n-m$ by applying $\textup{H}_{k-1}(C)$ to $\mathsf{D}_{k-1}$,\\[-0.cm]
  \begin{equation*}
  \begin{split}
  \delta_k \, \leq \, \bigg(\big(\ts\frac{c\nu_q^{3}}{\rho}\big)\,\Log(pn)^{\frac{5q+6}{2q}}\frac{1}{(n-k)\sqrt{k-1}}\bigg)
   \, + \, \bigg(\big( \ts\frac{cC\nu_q^7}{\rho^3}\big)(\Log(pn))^{r_q+\frac{3q+6}{2q}}\Log(k-1)\ts\frac{1}{(n-k)^{3/2}\sqrt{k-1}}\bigg).
  \end{split}
  \end{equation*}
  \normalsize
  Finally, to bound the sum $\delta_2+\cdots+\delta_{n-m}$, observe that
  \begin{equation*}
\sum_{k=2}^{n-m}\ts\frac{1}{(n-k)\sqrt{k-1}} \ \leq  \ \ts\frac{c\Log(n)}{n^{1/2}},
  \end{equation*}
  and
  \begin{equation*}
 \sum_{k=2}^{n-m}\ts\frac{1}{(n-k)^{3/2}\sqrt{k-1}} \ \leq  \ \ts\frac{c}{n^{1/2} m^{1/2}}.
  \end{equation*}
  Combining the last few steps leads to the stated result.\qed

\section{Proof of Theorem~\ref{thm:comp}}\label{sec:proofcomp}
Let $N$ be a positive integer that will be chosen later. Also, let $Z_1,\dots,Z_N$ be i.i.d.~copies of $Z$, and let $Y_1,\dots,Y_N$ be an independent sequence of i.i.d.~copies of $Y$. We will apply previous notations such as $S_{1:k}(X)$, $\delta_k^X(A)$, etc.~in a corresponding manner to the random vectors $Z_1,\dots,Z_N$, with $N$ playing the role of $n$ in the proof of Theorem~\ref{thm:main}. In particular, we have the following equalities in distribution for every choice of $N$,
\begin{equation*}
\begin{split}
S_{1:N}(Z) \ \, \scriptstyle\overset{\mathcal{L}}{\displaystyle=} \displaystyle\, \ Z,\\[0.2cm]
S_{1:N}(Y) \ \, \scriptstyle\overset{\mathcal{L}}{\displaystyle=}\displaystyle\, \ Y.
\end{split}
\end{equation*}
In order to re-use the proof of Theorem~\ref{thm:main}, the main part that needs to be revised is the moment matching argument in Section~\ref{sec:moment}. Specifically, the relation~\eqref{eqn:deltadiffmatch} must be modified, because in the current context, there is no guaranteed cancellation of the quadratic terms in the expansions~\eqref{eqn:deltaxexpand} and~\eqref{eqn:deltayexpand}. If we account for this detail in the reasoning leading up to~\eqref{eqn:deltadiffmatch}, then we have the following relation for every $k=1,\dots,N-1$,
\begin{equation}\label{eqn:newmatches}
\delta_k^Z(A)-\delta_k^Y(A) \ = \ \E[Q_k^Z(A)]-\E[Q_k^Y(A)] \ + \ \E[R_k^Z(A)]-\E[R_k^Y(A)].
\end{equation}

The remainder terms $\E[R_k^Z(A)]$ and $\E[R_k^Y(A)]$ can be handled in nearly the same manner as before in Section~\ref{sec:delta}. The only extra point to mention about bounding these remainders is the role of the parameters $\nu_2':=\max_{1\leq j\leq p}\|Z_{1j}\|_{\psi_2}$ and $\nu_2'':=\max_{1\leq j\leq p}\|Y_{1j}\|_{\psi_2}$. Under the reduction that $\Sigma^Y$ has all ones along the diagonal, Gaussianity implies that $\nu_2''$ is at most an absolute constant. Similarly, we have
$$(\nu_2')^2 \ \leq \ c \max_{1\leq j\leq p}\var(Z_{1j}) \ \leq \  c(1+\Delta^{\!\prime}),$$
and since we may assume without loss of generality that $\Delta'\leq 1$, it follows that $\nu_2'$ is also upper bounded by an absolute constant. 

To handle the difference of the quadratic terms $Q_k^Z(A)$ and $Q_k^Y(A)$ in~\eqref{eqn:newmatches}, observe that in the current context, the random hyperrectangle $\hat A_{k+1}$ defined previously in~\eqref{eqn:rhatdef} is independent of both $Z_k$ and $Y_k$. So, for $k=1,\dots,N-1$, we have
\begin{equation}
\begin{split}
\E\big[Q_k^Z(A)-Q_k^Y(A)\big]  & \ = 
\ \frac{1}{2}\,\E \bigg[\Big\langle \nabla^2 \varphi_{\e_k}(S_{1:k-1}(Z),\hat A_{k+1}), \ N^{-1} \big(Z_kZ_k\ttop -Y_kY_k\ttop\big) \Big\rangle\bigg]\\[0.2cm]
& \ = \ \frac{1}{2N}\Big\langle \E\Big[ \nabla^2 \varphi_{\e_k}(S_{1:k-1}(Z),\hat A_{k+1})\Big], \ \Sigma^Z-\Sigma^Y\Big\rangle.\label{eqn:gaussiansplit}
\end{split}
\end{equation}
Next, with regard to the Hessian of $\varphi_{\e_k}(\cdot,A)$, Lemma~\ref{lem:deriv} gives 
\begin{equation*}
\sup_{(s,A)\in\R^{^p}\times\mathscr{R}} \|\nabla^2 \varphi_{\e_k}(s,A)\|_1 \ \leq \ \ts\frac{c\Log(p)}{\e_k^2},
\end{equation*}
where $\e_k$ is defined in the current context as $\sqrt{\varrho(N-k)/N}$.
So, combining this with~\eqref{eqn:gaussiansplit} and H\"older's inequality, we have
\begin{equation}\label{eqn:holder2nd}
 \sup_{A\in\mathscr{R}} \Big|\E\big[Q_k^X(A)-Q_k^Y(A)\big] \Big| \ \leq \ \ts\frac{c \,\Log(p)\Delta^{\!\prime}}{\varrho(N-k)}.
\end{equation}

When re-using the proof of Theorem~\ref{thm:main}, the right side of~\eqref{eqn:holder2nd} should be added to the bounds on $\delta_k$ in the statements of Lemmas~\ref{lem:k2} and~\ref{lem:k1}. In addition, the induction hypothesis~\eqref{eqn:hyp} should be replaced with the inequality
 \begin{equation}\label{eqn:GcompHyp}
 \mathsf{D}_k \ \leq \ \frac{c\,\Log^4(pk)\Log(k)}{\varrho^{3/2} k^{1/2}} \ + \ \big(\ts\frac{c}{\varrho}\big)\Log(p)\Log(k)\Delta^{\!\prime}.                           
 \end{equation}
 Once these two updates are made, all of the corresponding steps in the proof of Theorem~\ref{thm:main} can be repeated to show there is an absolute constant $c>0$ such that the bound
\begin{equation}\label{eqn:almostGG}
\sup_{A\in\mathscr{R}}\Big| \P(Z\in A)-\P(Y\in A)\Big| \ \leq \  \displaystyle \frac{c\Log^4(pN)\Log(N)}{\varrho^{3/2} N^{1/2}} \ +  \big(\ts\frac{c}{\varrho}\big) \,\Log(p)\Log(N)\Delta^{\!\prime}
\end{equation}
holds for all $N$ and $p$.

The only remaining task is to choose $N$ in the bound~\eqref{eqn:almostGG}. To do this, first observe that  we may assume
\begin{equation}\label{eqn:simplify}
\ts\frac{1}{\varrho}\Log(p)\Delta^{\!\prime} \ \leq \ 1,
\end{equation}
for otherwise the theorem is true. Also, for purposes of simplification, note that there is an absolute constant $c>0$ such that
\begin{equation}\label{eqn:N13}
\ts \frac{\Log^4(pN)\Log(N)}{ N^{1/2}} \ \leq \frac{c\Log^4(p)}{N^{1/3}},
\end{equation}
with the exponent $1/3$ being unimportant.
If we choose $N$ such that
$$N^{1/3} \ = \ \lceil \ts\frac{1}{\sqrt \varrho \Delta^{\!\prime}}\Log^4(p)\rceil,$$
then~\eqref{eqn:almostGG} and~\eqref{eqn:N13} lead to
\begin{equation*}
\begin{split}
 \sup_{A\in\mathscr{R}}\Big| \P(Z\in A)-\P(Y\in A)\Big| & \ \leq \  \big(\ts\frac{c}{\varrho}\big) \,\Log(p)\Log(N)\Delta^{\!\prime}.
 \end{split}
\end{equation*}
Finally, observe that~\eqref{eqn:simplify} implies $\frac{1}{\sqrt{\varrho}\Delta^{\!\prime}}\log^4(p)\leq \varrho^{7/2}/(\Delta')^{5}\leq 1/(\Delta')^{5}$, and so the choice of $N$ implies $\Log(N)\leq c\Log(\frac{1}{\Delta^{\!\prime}})$, which leads to the stated result. \qed

\section{Proof of Theorem~\ref{thm:rho0}}
The argument is essentially the same as that of the proof of Theorem~\ref{thm:comp} up to the bound~\eqref{eqn:almostGG}, with correspondences of notation given by $n\leftrightarrow N$, $S_{1:n}(X)\leftrightarrow S_{1:N}(Z)$, and $\Delta^{\!\prime}\leftrightarrow \Delta$. 
 One slight difference to note is that the induction hypothesis~\eqref{eqn:GcompHyp} should be modified to account  for the parameters $\omega_q$ and $r_q$ as
 $$\mathsf{D}_k \ \leq \ \frac{c\,\omega_q^4\,\Log^{r_q}(pk)\Log(k)}{\varrho^{3/2} k^{1/2}} \ + \ \big(\ts\frac{c}{\varrho}\big)\Log(p)\Log(k)\Delta. $$
\qed

\section{Proof of Theorem~\ref{thm:boot}}\label{sec:bootproof} By Theorem~\ref{thm:rho0}, it suffices to bound the distance between the distribution of $Y$ and the conditional distribution of $\ts\frac{1}{\sqrt n}\sum_{i=1}^n X_i^{\star}$ given $X_1,\dots,X_n$. 
Since the sum $\ts\frac{1}{\sqrt n}\sum_{i=1}^n X_i^{\star}$ is conditionally Gaussian with mean zero and covariance matrix $\hat\Sigma$, it follows from Theorem~\ref{thm:comp} that
\begin{equation}\label{eqn:thm:boot:proof}
\sup_{A\in\mathscr{R}}\Bigg|\P\Big(\ts\frac{1}{\sqrt n}\sum_{i=1}^n X_i^{\star}\in A\Big|X\!\Big)-\P\big(Y\in A\big)\Bigg|
 \ \leq \  
  \displaystyle \big(\ts\frac{c}{\varrho}\big)\Log(p)\Log\big(\ts\frac{1}{\hat\Delta}\big) \hat\Delta,
\end{equation}
where we define the random variable $\hat\Delta=\|\hat\Sigma-\Sigma^Y\|_{\infty}$. To bound $\hat\Delta$, first note that we may assume $n^{1/2}\geq \log^{r_q}(pn)\log(n)$ without loss of generality, for otherwise the stated result is true. In particular, this implies $n^{1/2}\geq \log^{2/q}(pn)\log^{q/2}(n)$, and so it follows from part (v) of Lemma~\ref{lem:tail} that the event
\begin{equation}\label{eqn:splitdeltahat}
\hat\Delta \ \leq \ c \omega_q^2\sqrt{\ts\frac{\log(pn)}{n}} + \Delta
\end{equation}
holds with probability at least $1-c/n$.

Next, in order to deal with the expression $\log(\frac{1}{\hat\Delta})\hat\Delta$, we may assume that the right side of~\eqref{eqn:splitdeltahat} is at most $1/3$ without loss of generality, for again, the stated result would be true otherwise. Also, note that that the function $z\mapsto \log(\frac{1}{z}) z$ is increasing on the interval $[0,1/3]$, with the expression $\log(\frac{1}{0})0$ understood as 0. Consequently, if $x>0$ and $y,z\geq 0$ are numbers satisfying $z\leq x+y\leq 1/3$, it follows that
$$\log(\ts\frac{1}{z}) z \ \leq \ \log(\ts\frac{1}{x})(x+y).$$
So, by viewing~\eqref{eqn:splitdeltahat} in the form $z\leq x+y$, and recalling the bound $\omega_q\geq \sqrt{1/2}$ from page \pageref{page:prelim}, it follows that the event
\begin{equation}\label{eqn:twopartbound}
\log\big(\ts\frac{1}{\hat\Delta}\big)\hat\Delta \ \leq \ c\log(n)\Big(c \omega_q^2\sqrt{\ts\frac{\log(pn)}{n}}  \ + \ \Delta\Big)
\end{equation}
holds with probability at least $1-c/n$.
  Lastly, note that when the right side of~\eqref{eqn:twopartbound} is multiplied by $(1/\varrho)\log(p)$, the resulting value is of at most the same order as the error arising from the Gaussian approximation result in Theorem~\ref{thm:rho0}. This completes the proof.
 \qed

 \section{Proof of Theorem~\ref{thm:bootemp}}\label{sec:bootempproof} As in the proof of the multiplier bootstrap approximation result, it is enough to bound the distance between the distribution of $Y$ and the conditional distribution of  $\ts\frac{1}{\sqrt n}\sum_{i=1}^n (X_i^*-\bar X)$ given $X_1,\dots,X_n$. For a generic random variable $V$, define the conditional Orlicz norm $\|V\|_{\psi_q|X}=\inf\{t>0\, |\, \E[\exp(|V|^q/t^q)|X]\leq 2\}$, and define the bootstrap counterpart of $\omega_q$ as
\begin{equation}
\hat \omega_q = \max_{1\leq j\leq p} \big\|X_{1j}^*-\bar X_j\big\|_{\psi_q|X}.
\end{equation}
Also, recall the quantity $\hat\Delta=\|\hat\Sigma-\Sigma^Y\|_{\infty}$.
Then, Theorem~\ref{thm:rho0} implies that the bound
\small
\begin{equation}
\sup_{A\in\mathscr{R}}\Big|\P\Big(\ts\frac{1}{\sqrt n}\sum_{i=1}^n (X_i^*-\bar X) \in A\Big|X\!\Big)-\P\big(Y\in A\big)\Big|
 \   \ \leq \ \displaystyle\frac{ c\, \hat \omega_q^4 \,\Log^{r_q}(pn)\,\Log(n)}{\varrho^{3/2}\, n^{1/2}} \ +
  \  \displaystyle \big(\ts\frac{c}{\varrho}\big)\Log(p)\Log(n) \hat \Delta
\end{equation}
\normalsize
holds almost surely.  To bound the quantity $\hat\omega_q$, it is straightforward to check first that
\begin{equation}
\hat\omega_q \ \leq \ c\max_{1\leq i\leq n}\|X_i\|_{\infty},
\end{equation}
and so by part (ii) Lemma~\ref{lem:tail}, it follows that the bound
\begin{equation}
\hat\omega_q \ \leq \ c\omega_q \Log^{1/q}(pn)
\end{equation}
holds with probability at least $1-c/n$. The proof is completed by combining this bound on $\hat\omega_q$ with the bound on $\hat\Delta$ given in~\eqref{eqn:splitdeltahat}.\qed

 \section{Background results}\label{sec:background}

The facts about random vectors in the following lemma are mostly standard. When using these facts, it may be helpful to note that the conditions of Theorem~\ref{thm:main} are a special case of those in Theorem~\ref{thm:rho0}, and also, that the equality $\omega_q=\nu_q$ holds under the conditions of Theorem~\ref{thm:main}.

  \begin{lemma}\label{lem:tail} 
There is an absolute constant $c>0$ such that the following statements hold for all $n,p\geq 1$ and $q\in\{1,2\}$, provided that the conditions of Theorem~\ref{thm:rho0} hold, and $\var(Y_{j})=1$ for all $1\leq j\leq p$:
\begin{enumerate}[(i)]
\item The expectation of $\|X_1\|_{\infty}^3$ satisfies
 \begin{equation*}
\E[\|X_1\|_{\infty}^3] \ \leq \ c\omega_q^3 \Log^{3/q}(p).
\end{equation*}
\item The event
$$\max_{1\leq i\leq n}\|X_i\|_{\infty} \ \leq \ c\omega_q\log^{1/q}(pn),$$
holds with probability at least $1-c/n$.\\[-0.2cm]
\item If $t=c(\omega_q+1) \Log^{1/q}(pn)n^{-1/2}$, then
\begin{equation*}
\P(\ts\frac{1}{\sqrt n}\|X_1-Y_1\|_{\infty}\geq t) \ \leq \ \ts\frac{c}{n}. 
\end{equation*}
\item If $1\leq m\leq n$, and the vectors $\zeta=S_{1:n}(X)$ and $\xi=S_{1:n-m}(X)+S_{n-m+1:n}(Y)$ are as in the proof of Lemma~\ref{lem:start}, then the following bound holds when $t'=c(\omega_q+1) \sqrt{\frac{m}{n}}\Log^{1/q}(pn)$,
\begin{equation}\label{eqn:2ndtail}
\P(\|\zeta-\xi\|_{\infty}\geq t') \ \leq \ \ts\frac{c}{n}. 
\end{equation}
\item If the matrix $\hat\Sigma$ is as defined in~\eqref{eqn:sigmahatdef}, then the event
$$\|\hat\Sigma-\E[X_1X_1\ttop]\|_{\infty} \ \leq \ c \omega_q^2\Big(\sqrt{\ts\frac{\log(pn)}{n}}+\ts\frac{\log^{2/q}(n)\log^{2/q}(pn)}{n}\Big)$$
holds with probability at least $1-c/n$.
\end{enumerate}
 \end{lemma}
 \proof Items (i), (ii), and (iii) can be proved using the results in~\citep[][Section 2.2]{Wellner} and~\citep[][Sections 2.5-2.6]{Vershynin}. Item (v) follows from Theorem 4.2 of~\citep{KC:2018}.\qed \\

The next result is known as Nazarov's Gaussian anti-concentration inequality, which originates from the paper~\citep{Nazarov:2003}, and was further elucidated in~\citep[][Theorem 1]{CCK:Nazarov}. For this lemma and the next one, recall that the sets $A^t$, $A^{-t}$ and $\partial A(t)$ associated with a hyperrectangle $A$ are defined in the notation paragraph of Section~\ref{sec:intro}.

   \begin{lemma}\label{lem:anti}
There is an absolute constant $c>0$ such that the following holds for all $p$: Let $\xi\in\R^p$ be a Gaussian random vector, and suppose that  $\varsigma=\min_{1\leq j\leq p}\sqrt{\var(\xi_j)}$ is positive. Then, for any $t>0$,
\begin{equation*}
\sup_{A\in\mathscr{R}}\P(\xi\in \partial A(t)) \ \leq \ \ts\frac{ct}{\varsigma}\sqrt{\Log(p)}.
\end{equation*}
   \end{lemma}  
\noindent As a clarification on this version of Nazarov's inequality, it is worth noting that the result is more commonly stated in terms of probabilities of the form $\P(\xi\in \partial C_b(t))$, where we define the ``corner set''
\smash{$C_b=
\{x\in\R^p|\, x_j\leq b_j \text{ for all } j=1,\dots,p\}$}
associated with a fixed vector $b=(b_1,\dots,b_p)$. The lemma above can be reduced to the more common version by noting that if $A=\prod_{j=1}^p [a_j,b_j]$ and $a=(a_1,\dots,a_p)$, then we have the inclusion
\begin{equation}\label{eqn:subsetbound}
\{\xi\in \partial A(t)\} \ \subset \ \Big(\big\{-\xi\in \partial C_{-a}(t)\big\}\cup \big\{\xi\in\partial C_b(t)\big\}\Big).
\end{equation}
 Thus, the more common version of Nazarov's inequality can be used after a union bound is applied to~\eqref{eqn:subsetbound}.

The next result is often used for scalar random variables, but it seems to be stated less frequently in the case of random vectors. Also, this result remains valid when the boundary set $\partial A(t)$ is defined with respect to an arbitrary norm $\|\cdot\|$ on $\R^p$ (rather than just the $\|\cdot\|_{\infty}$ norm used in the definition given in Section~\ref{sec:intro}).

\begin{lemma}\label{lem:coupling}
Let $\|\cdot\|$ be any norm on $\R^p$, and let $\zeta,\xi\in\R^p$ be any two random vectors. Then, the following inequality holds for any Borel set $A\subset\R^p$, and any $t>0$,
\[
|\P(\zeta \in A)-\P(\xi\in A )| \ \leq \  \P\big(\xi\in \partial A(t)\big) \ + \   \P\big(\|\zeta-\xi\|\geq t\big).
\]
\end{lemma}

\proof Let $\delta=\zeta-\xi$ and observe that
$$  \P\big(\xi \in A^{-\|\delta\|}\big) \  \leq \  \P(\zeta \in A) \ \leq \   \P\big(\xi \in A^{\|\delta\|}\big).
$$
This implies
 \begin{equation*}
 \begin{split}
  \left|\P(\zeta\in A)-\P(\xi\in A) \right|  & \ \leq \ \P\Big( \xi \in  \big(A^{\|\delta\|}\backslash  A^{-\|\delta\|}\big)\Big) 
   \ \leq  \  \P\Big( \xi \in  (A^{t}\backslash  A^{-t})\Big) \ + \ \P(\|\delta\|\geq t).
\end{split}
\end{equation*}
\qed

For the following lemma, recall that $\hat q_{1-\alpha}$ denotes the $1-\alpha$ quantile of the distribution function $t\mapsto \P(M_n^*\leq t|X)$ defined on page~\pageref{page:quantile}.

\begin{lemma}\label{lem:quantiles}
There is an absolute constant $c>0$ such that the following holds for all $n,p\geq 1$, $q\in\{1,2\}$, and $\alpha\in(0,1)$:
If the random vectors $X_1,\dots,X_n\in\R^p$ satisfy the conditions of Theorem~\ref{thm:main}, then
\begin{equation}\label{eqn:quantiles}
\big|\P(M_n \leq \hat q_{1-\alpha})-(1-\alpha)\big| \ \leq \  \frac{ c\,\nu_q^{4}\,\Log^{s_q}(pn)\,\Log(n)}{\rho^{3/2}\, n^{1/2}}.
\end{equation}
\end{lemma}

\proof  
Let the Gaussian random vector $Y\in\R^p$ be as in the statement of Theorem~\ref{thm:main}, and let $q_{1-\alpha}^Y$ denote the $1-\alpha$ quantile of the random variable $\|Y\|_{\infty}$. 
 In addition, let $\eta$ denote a number of the form
$$\eta  \ =  \ \frac{c\nu_q^{4}\,\Log^{s_q}(pn)\,\Log(n)}{\rho^{3/2}n^{1/2}},$$
where $c>0$ is an absolute constant that is at least as large as $C$ in Theorem~\ref{thm:main}, and at least as large as any of the instances of $c$ in the statements of the other theorems of Section~\ref{sec:main}. 
First, we will show that 
 \begin{equation}\label{eqn:claim1}
  \P(M_n\leq \hat q_{1-\alpha}) \ \geq  \ 1-\alpha-4\eta-c/n.
 \end{equation}
In this portion of the proof, we may assume that $1-\alpha-4\eta>0$, since the inequality above is  true otherwise. For this reason, it will make sense to work with quantiles corresponding to $1-\alpha-4\eta$.
 Due to the Gaussianity of $Y$, the distribution of $\|Y\|_{\infty}$ is absolutely continuous with respect to Lebesgue measure, and so \smash{$\P(\|Y\|_{\infty}\leq q_{1-\alpha-3\eta}^Y)=1-\alpha-3\eta$.} Accordingly, the Gaussian approximation result of Theorem~\ref{thm:main} gives
\begin{equation}\label{eqn:Gsandwich}
1-\alpha-4\eta\  \ \leq \ \P(M_n\leq q_{1-\alpha-3\eta}^Y) \ \leq \ 1-\alpha-2\eta. 
\end{equation}
The first inequality above implies $q_{1-\alpha-4\eta}\leq q_{1-\alpha-3\eta}^Y$, and then combining this with the second inequality above gives
\begin{equation}\label{eqn:Gstep}
\begin{split}
\P(M_n\leq q_{1-\alpha-4\eta}) 
&  \ \leq \  1-\alpha-2\eta.
\end{split}
\end{equation}
Next, we may use this bound and the bootstrap approximation result in Theorem~\ref{thm:bootemp} to conclude that the event
\begin{equation}
\begin{split}
\P(M_n^*\leq q_{1-\alpha-4\eta}|X) 
& \ \leq \ 1-\alpha-\eta
\end{split}
\end{equation}
holds with probability at least $1-c/n$. Consequently, the event $G=\{\hat q_{1-\alpha}> q_{1-\alpha-4\eta}\}$ must occur with at least the same probability.
Next, observe
\begin{equation}
\begin{split}
\P(M_n> \hat q_{1-\alpha}) & \ \leq \  \P\Big(\{M_n> \hat q_{1-\alpha}\}\cap G\Big) + \ts\frac{c}{n}\\
&  \ \leq \ \P(M_n> q_{1-\alpha-4\eta})+\ts\frac{c}{n}\\
& \ \leq \ \alpha+4\eta+\ts\frac{c}{n}.
\end{split}
\end{equation}
Rearranging this implies the claim~\eqref{eqn:claim1}. To complete the proof, a similar argument can be used to show that
\begin{equation}
\P(M_n\leq \hat q_{1-\alpha})\leq 1-\alpha+3\eta+c/n.
\end{equation}
In essence, this involves showing that the event $G'=\{\hat q_{1-\alpha}\leq q_{1-\alpha+\eta}\}$ holds with probability at least $1-c/n$, and then using $\P(M_n\leq \hat q_{1-\alpha})  \ \leq \  \P(\{M_n\leq \hat q_{1-\alpha}\}\cap G') + \frac{c}{n}$. \qed

~\\

The last lemma in this section handles some low-level details in the proof of Lemma~\ref{lem:k2}.

\begin{lemma}\label{lem:bkbound}
There is an absolute constant $c>0$ such that the following holds for all $n,p\geq 1$. If $\e>0$ and $\ve=4\e\sqrt{\Log(pn)}$, then the inequality below holds for any pair \smash{$(s,A)\in \R^p\times\mathscr{R}$} such that $s\not\in \partial A(\ve)$,
\begin{equation*}
\big\|\nabla^3\varphi_{\e}(s,A)\big\|_1\leq \ts\frac{c}{\e^3 pn}.
\end{equation*}
\end{lemma}
\proof Observe that
\begin{equation}\label{eqn:fullderiv}
\big\|\nabla^3\varphi_{\e}(s,A)\big\|_1 = \sum_{j_1,j_2,j_3}
|\partial_{j_1,j_2,j_3}\varphi_{\e}(s,A)|
\end{equation}
where $j_1$, $j_2$, and $j_3$ may take any values in $\{1,\dots,p\}$,  and we use the shorthand notation $\partial_{j_1,j_2,j_3}=  \ts\frac{\partial^3}{\partial s_{j_1}\partial s_{j_2}\partial s_{j_3}}$. To work with this sum, we will focus on the case when $A$ can be represented as a  product of compact intervals $A=\prod_{j=1}^p [a_j,b_j]$, since the  argument is virtually the same for other types of hyperrectangles.

There are three possible forms that the summands in~\eqref{eqn:fullderiv} can take. First, when $j_1$, $j_2$, and $j_3$ are all distinct, it follows from the formula~\eqref{eqn:phiprod2} that 
\begin{equation*}
\small
\begin{split}
\big|\partial_{j_1,j_2,j_3}\varphi_{\e}(s,A)\big| & \  = \ 
 \frac{1}{\e^3} \displaystyle\prod_{i=1}^3 \Big|\phi\big(\ts\frac{b_{j_i}-s_{j_i}}{\e}\big)-\phi\big(\ts\frac{a_{j_i}-s_{j_i}}{\e}\big)\Big| \displaystyle \prod_{l\not\in\{j_1,j_2,j_3\}} \Big(\Phi(\ts\frac{b_{l}-s_{l}}{\e})-\Phi(\ts\frac{a_{l}-s_{l}}{\e})\Big).
 \end{split}
\end{equation*}
Second, when exactly two of the indices are the same (say $j_1\neq j_2$ and $j_1=j_3$), we have
\begin{equation*}
\small
\begin{split}
& \big|\partial_{j_1,j_2,j_1}\varphi_{\e}(s,A)\big| \ = \  \\[0.2cm]
& \frac{1}{\e^3}\Big|\phi'
\big(\ts\frac{b_{j_1}-s_{j_1}}{\e}\big)-\phi'\big(\ts\frac{a_{j_1}-s_{j_1}}{\e}\big)\Big|\Big|\phi\big(\ts\frac{b_{j_2}-s_{j_2}}{\e}\big)-\phi\big(\ts\frac{a_{j_2}-s_{j_2}}{\e}\big)\Big|
  \displaystyle \prod_{l\not\in\{j_1,j_2\}} \Big(\Phi(\ts\frac{b_{l}-s_{l}}{\e})-\Phi(\ts\frac{a_{l}-s_{l}}{\e})\Big).
 \end{split}
\end{equation*}
Third, when $j_1=j_2=j_3$, we have
\begin{equation*}
\small
%\begin{split}
 \big|\partial_{j_1,j_1,j_1}\varphi_{\e}(s,A)\big|  \ = \ 
 \frac{1}{\e^3}\Big|\phi''
\big(\ts\frac{b_{j_1}-s_{j_1}}{\e}\big)-\phi''\big(\ts\frac{a_{j_1}-s_{j_1}}{\e}\big)\Big|
  \displaystyle \prod_{l\neq j_1} \Big(\Phi(\ts\frac{b_{l}-s_{l}}{\e})-\Phi(\ts\frac{a_{l}-s_{l}}{\e})\Big).
% \end{split}
\end{equation*}
For the subsequent discussion, it will be helpful to note that the non-negative functions $\Phi$, $\phi$, $|\phi'|$, and $|\phi''|$ are all uniformly bounded over $\R$ by an absolute constant.

To proceed, consider the two ways in which the condition  $s\not\in\partial A(\ve)$ can hold. Namely, $s$ must either (i) lie in $A\setminus \partial A(\ve)$, or (ii) lie in  $\R^p\setminus A^{\ve}$. The first case implies in particular that the inequalities $s_{j_1}-a_{j_1}\geq \ve$ and $b_{j_1}-s_{j_1}\geq \ve$ hold for any choice of $j_1$. In addition, note that the inequality $\max\{\phi(x),|\phi'(x)|,|\phi''(x)|\}\leq x^2\phi(x)$ holds for all $|x|\geq 1$. Consequently, in case (i), for any of the triples $(j_1,j_2,j_3)$, it follows that if $\ve=4\e\sqrt{\Log(pn)}$, then
\begin{equation}
\begin{split}
\big|\partial_{j_1,j_2,j_3}\varphi_{\e}(s,A)\big|  & \ \leq  \ts\frac{c\Log(pn)}{\e^3}\phi\Big(4\sqrt{\Log(pn)}\Big)\\[0.2cm]
& \ \leq \ \ts\frac{c\Log(pn)}{\e^3}\frac{1}{(pn)^{8}}.\label{eqn:last3derivbound}
\end{split}
\end{equation}
Moreover, since there are $p^3$ terms in the sum~\eqref{eqn:fullderiv}, we have $\|\nabla^3\varphi_{\e}(s,r)\|_1\leq \frac{c}{\e^3 (pn)^4}.$ 

Finally, consider case (ii) where $s\in \R^p\setminus A^{\ve}$. There must be at least one $j\in\{1,\dots,p\}$ such that $s_j\not\in [a_j-\ve,b_j+\ve]$. If this $j$ belongs to the triple $(j_1,j_2,j_3)$, then the bound~\eqref{eqn:last3derivbound} holds. Alternatively, if $j$ does not belong to this triple, then the identity $1-\Phi(x)=\Phi(-x)$ holding for all $x\in\R$, and the bound 
%\begin{equation}
$\Phi(-x) \ \leq \ e^{-x^2/2}$
%\end{equation}
holding for all $x\geq 1$,
implies
$$\Phi(\ts\frac{b_{j}-s_{j}}{\e})-\Phi(\ts\frac{a_{j}-s_{j}}{\e}) \ \leq  \  \ts\frac{c}{(pn)^{8}}.$$
Altogether, we conclude that
$\|\nabla^3\varphi_{\e}(s,A)\|_1\leq \frac{c}{\e^3 (pn)^4}$.\qed

%%  MR numbers will be added by VTeX.                      %%
%%                                                         %%
%%  Use \cite{...} to cite references in text.             %%
%%                                                         %%
%%%%%%%%%%%%%%%%%%%%%%%%%%%%%%%%%%%%%%%%%%%%%%%%%%%%%%%%%%%%%

%% if your bibliography is in bibtex format, uncomment commands:
%\bibliographystyle{imsart-number} % Style BST file (imsart-number.bst or imsart-nameyear.bst)
%\bibliography{bibliography}       % Bibliography file (usually '*.bib')

%% or include bibliography directly:
\references

\end{document}